\documentclass[a4paper,10pt]{article}

\usepackage{preambule}
\usepackage{chngcntr}
\theoremstyle{thmstyletwo}%

\numberwithin{equation}{section}

\newcommand{\eg}{\textit{e.g.}~}

\newcommand{\B}{\mathbf{B}}
\newcommand{\CBold}{\mathbf{C}}
\newcommand{\D}{\mathbf{D}}
\newcommand{\G}{\mathbf{G}}
\renewcommand{\H}{\mathbf{H}}
\newcommand{\J}{\mathbf{J}}
\renewcommand{\K}{\mathbf{K}}
\newcommand{\LBold}{\mathbf{L}}
\newcommand{\NBold}{\mathbf{N}}
\newcommand{\MBold}{\mathbf{M}}
\newcommand{\T}{\mathbf{T}}
\newcommand{\U}{\mathbf{U}}
\newcommand{\V}{\mathbf{V}}
\newcommand{\W}{\mathbf{W}}
\newcommand{\X}{\mathbf{X}}
\newcommand{\Y}{\mathbf{Y}}

\newcommand{\ZBold}{\mathbf{Z}}
\newcommand{\SBold}{\mathbf{S}}
\newcommand{\RBold}{\mathbf{R}}
\newcommand{\QBold}{\mathbf{Q}}

\newcommand{\SigmaBold}{\bm{\Sigma}}
\newcommand{\OmegaBold}{\bm{\Omega}}

\newcommand{\PhiBold}{\bm{\Phi}}

\RequirePackage[normalem]{ulem} 
\RequirePackage{color}\definecolor{RED}{rgb}{1,0,0}\definecolor{BLUE}{rgb}{0,0,1} 
\RequirePackage{listings} 
\RequirePackage{color} 
\lstdefinelanguage{DIFcode}{ 
  moredelim=[il][\color{red}\sout]{\%DIF\ <\ }, 
  moredelim=[il][\color{blue}\uwave]{\%DIF\ >\ } 
} 
\lstdefinestyle{DIFverbatimstyle}{ 
	language=DIFcode, 
	basicstyle=\ttfamily, 
	columns=fullflexible, 
	keepspaces=true 
} 
\lstnewenvironment{DIFverbatim}{\lstset{style=DIFverbatimstyle}}{} 
\lstnewenvironment{DIFverbatim*}{\lstset{style=DIFverbatimstyle,showspaces=true}}{} 

\graphicspath{{Fig/}}

\newcommand{\Proj}{\Pi}

\title{A unified error analysis for randomized low-rank approximation with application to data assimilation}

\author{
    Alexandre Scotto Di Perrotolo
    \thanks{IRT Saint Exupéry, 3 Rue Tarfaya, 31400 Toulouse, France (\href{mailto:alexandre.scotto@irt-saintexupery.com}{alexandre.scotto@irt-saintexupery.com}).} \and
    Youssef Diouane
    \thanks{GERAD and Department of Mathematics and Industrial Engineering, Polytechnique Montr\'eal, Canada (\href{mailto:youssef.diouane@polymtl.ca}{youssef.diouane@polymtl.ca}).} \and
    Selime G\"urol
    \thanks{CERFACS, 42 Avenue Gaspard Coriolis, F-31057 Toulouse Cedex 01, France (\href{mailto:selime.gurol@cerfacs.fr}{selime.gurol@cerfacs.fr}).} \and
    Xavier Vasseur\thanks{4 rue Jean-Pierre Petit, F-31700 Blagnac, France (\href{mailto:xavier.vasseur@gmail.com}{xavier.vasseur@gmail.com}).}
}

\begin{document}

\maketitle

\begin{abstract}
    Randomized algorithms have proven to perform well on a large class of numerical linear algebra problems. Their theoretical analysis is critical to provide guarantees on their behaviour, and in this sense, the stochastic analysis of the randomized low-rank approximation error 
    plays a central role. Indeed, several randomized methods for the approximation of dominant eigen- or singular modes can be rewritten as low-rank approximation methods. However, despite the large variety of algorithms, the existing theoretical frameworks for their analysis rely on a specific structure for the covariance matrix that is not adapted to all the algorithms. We propose a general framework for the stochastic analysis of the low-rank approximation error in Frobenius norm for centered and non-standard Gaussian matrices. Under minimal assumptions on the covariance matrix, we derive accurate bounds both in expectation and probability. Our bounds have clear interpretations that enable us to derive properties and motivate practical choices for the covariance matrix resulting in efficient low-rank approximation algorithms. The most commonly used bounds in the literature have been demonstrated as a specific instance of the bounds proposed here, with the additional contribution of being tighter. Numerical experiments related to data assimilation further illustrate that exploiting the problem structure to select the covariance matrix improves the performance as suggested by our bounds.
\end{abstract}

\medskip


\medskip

\textbf{MSCcodes:} \textbf{65F55, 65F99, 65K99, 68W20.}

\section{Introduction}
    Let $\A \in \Rmn$ be an arbitrary real matrix and consider the following minimization problem,
    \begin{equation} \label{eq:low_rank_prob}
        \min_{\mathbf{B} \in \Rmn} \norm{\A - \mathbf{B}}_F \quad \text{subject to} \quad \rank(\mathbf{B}) = k \leq \min\set{m, n},
    \end{equation}
    where $\norm{\cdot}_F$ denotes the Frobenius norm. This problem, referred to as the \textit{low-rank approximation} problem, is a key ingredient in numerous applications in data analysis and scientific computing including principal component analysis~\cite{RokhlinSzlamEtAl_RandomizedAlgorithmPrincipal_2010}, data compression~\cite{Mahoney_RandomizedAlgorithmsMatrices_2010} and approximation algorithms for partial differential and integral equations~\cite{Hackbusch_HierarchicalMatricesAlgorithms_2015}, to name a few. Its solution~\cite{EckartYoung_ApproximationOneMatrix_1936} is obtained from the order $k$ truncated singular value decomposition (SVD) of $\A$.

    Let $\U \in \Rmm$ and $\V \in \Rnn$ be the real orthogonal matrices containing the left and right singular vectors of $\A$ respectively, and $\SigmaBold \in \Rmn$ the matrix whose diagonal contains the singular values of $\A$ with the convention $\sigma_1 \geq \dots \geq \sigma_{\min\set{m, n}}$. For a given integer $k \leq \min\set{m, n}$, we consider the following partitioning of the SVD of $\A$,
    \begin{equation} \label{eq:partitioning_svd}
        \A =
        \begin{bmatrix}
            \U_k & \overline{\U}_k
        \end{bmatrix}
        \begin{bmatrix}
            \SigmaBold_k & \\
            & \overline{\SigmaBold}_k
        \end{bmatrix}
        \begin{bmatrix}
            \V_k & \overline{\V}_k
        \end{bmatrix}\tra,
    \end{equation}
    where $\U_k \in \Rmk$, $\overline{\U}_k \in \R^{m \times (m-k)}$, $\V_k \in \Rnk$, $\overline{\V}_k \in \R^{n \times (n-k)}$, $\SigmaBold_k \in \Rkk$ and $\overline{\SigmaBold}_k \in \R^{(m-k) \times (n-k)}$. Defining $\A_k = \U_k \SigmaBold_k \V_k\tra$ and $\overline{\A}_k = \overline{\U}_k \overline{\SigmaBold}_k \overline{\V}_k\tra$, one can rewrite $\A = \A_k + \overline{\A}_k$ where the matrix $\A_k$ is precisely the minimizer of~\eqref{eq:low_rank_prob} with the corresponding optimal value being $\norm{\A - \A_k}_F = \norm{\overline{\SigmaBold}_k}_F$.

    The matrix $\A_k$ can be computed with any SVD algorithm enabling a truncation mechanism~\cite[Section 9.6]{GolubVanLoan_MatrixComputations_2013}. However, for large-scale problems, the classical approaches become prohibitively expensive or not even inapplicable if $\A$ is not stored explicitly as a matrix. In these situations, one is rather interested in computing $\widehat{\A}_k$, a reasonably accurate approximation of $\A_k$, which is typically expected to be optimal up to a small factor $\varepsilon > 0$, that is $\widehat{\A}_k$ satisfies,
    \begin{equation*}
        \norm{\A - \widehat{\A}_k}_F \leq
        (1 + \varepsilon) \, 
        \norm{\A - \A_k}_F \displaystyle =   (1 + \varepsilon) \, \norm{\overline{\SigmaBold}_k}_F.
    \end{equation*}
    Using $\A_k = \U_k \U_k\tra \A = \Proj(\U_k) \A$, a natural strategy consists of searching for a low-rank approximation of $\A$ of the form $\widehat{\A}_k = {\Proj}(\ZBold) \A$. If $\range(\cdot)$ denotes the column space, then the matrix $\ZBold$ should be computed such that $\range(\ZBold)$ accurately approximates $\range(\U_k)$, and algorithms computing such matrices are referred to as \textit{range-finders}. In this case, the low rank approximation error takes the particular form
    \begin{equation} \label{eq:metric}
        \norm{[\I_m -  \Proj(\ZBold)] \A}_F.
    \end{equation}
    
Randomized algorithms for computing $\ZBold$ have been proposed. Following~\cite{WoolfeLibertyEtAl_FastRandomizedAlgorithm_2008}, the authors in~\cite{HalkoMartinssonEtAl_FindingStructureRandomness_2011} proposed an efficient algorithm now widely known as the Randomized SVD (RSVD) algorithm, which considers $\ZBold = \A \G$, where $\G$ is a matrix whose columns are independently sampled from an $n$-variate standard Gaussian distribution. More importantly, they provided expectation and probability bounds for \eqref{eq:metric} which were critical to identify key properties of the RSVD algorithm, and give confidence in its performance.
    In a subsequent study, the authors in~\cite{Gu_SubspaceIterationRandomization_2015} extended the analysis to the full randomized subspace iteration~\cite{RokhlinSzlamEtAl_RandomizedAlgorithmPrincipal_2010} and provided stochastic error bounds when $\ZBold = \A (\A\tra\A)^q \G$ and $q \geq 0$.
    
    Beyond the RSVD, various randomized algorithms were proposed to address more elaborated eigenvalue and singular value problems.     In~\cite{SaibabaLeeEtAl_RandomizedAlgorithmsGeneralized_2016}, the authors proposed methods to compute dominant eigenpairs of an Hermitian matrix $\A_\text{H}$ with respect to an Hermitian positive definite matrix $\B$. This problem is actually equivalent to finding the low rank approximation of $\A_\text{B} \equiv \B\halfinv \A_\text{H} \B\halfinv$, and several randomized algorithms were proposed, constructing the low rank approximation from matrices of the form $\ZBold = \A_\text{B}^q \B\half \G$ where $q\in \{1,2\}$. For $q=1$, an expectation bound  has been proposed in~\cite[Theorem 1]{SaibabaLeeEtAl_RandomizedAlgorithmsGeneralized_2016}  

but no probability bounds.
    In~\cite{SaibabaHartEtAl_RandomizedAlgorithmsGeneralized_2021}, the authors derived algorithms for computing truncated generalized singular value decomposition~\cite{VanLoan_GeneralizingSingularValue_1976}. The proposed decomposition involves two symmetric positive-definite matrices $\SBold \in \Rmm$ and $\T \in \Rnn$ and reduces to the computation of a low rank approximation of $\A_\text{S,T} \equiv \SBold^{1/2} \A \T^{-1/2}$. The authors provide a randomized algorithm based on the power iteration scheme, that is $\ZBold = (\A_\text{S,T} \, \A\tra_\text{S,T})^q \A_\text{S,T} \T^{1/2} \G$ and $q \geq 0$. The derived probability bounds show that the approximation error vanishes asymptotically with $q$. However, the bounds predict a deviation from the optimum error that depends on the dimension $n$ and the 2-norm conditioning of $\T$, which seems to severely overestimate the actual error observed in the experiments. Also no bounds in expectation have been derived.

Note that all the studied problems reduce to construct a low rank approximation of a given matrix $\A$ out of a randomly drawn matrix of the form $\ZBold = (\A\A)^q\A \LBold \G$, where $\LBold \in \Rnn$ depends on the specific problem considered. This suggests that providing an accurate stochastic analysis of \eqref{eq:metric} for matrices $\ZBold$ of this form is highly relevant and has not yet been thoroughly studied in the literature.
    Recently, a generalization of the RSVD to the infinite-dimensional case, where $\A$ is replaced by a linear differential operator, has been proposed in~\cite{BoulleTownsend_LearningEllipticPartial_2022}. From their theoretical analysis in the infinite-dimensional setting, they were able to propose error bounds~\cite{BoulleTownsend_GeneralizationRandomizedSingular_2022} for~\eqref{eq:metric} in the more general setting $\ZBold = \A \OmegaBold$ where the columns of $\OmegaBold$ are independently sampled from an $n$-variate centered Gaussian distribution with covariance matrix $\CBold \in \Rnn$. Using standard properties of the multivariate Gaussian distribution, this case is equivalent to $\ZBold = \A \LBold \G$, where $\LBold \in \Rnn$ is such that $\CBold = \LBold \LBold\tra$. Numerically, their approach showed that a relevant choice for $\CBold$ could improve the performance. However, plugging $\LBold = \I_n$ in~\cite[Theorem 2]{BoulleTownsend_GeneralizationRandomizedSingular_2022} yields poorer bounds compared to the ones in~\cite{HalkoMartinssonEtAl_FindingStructureRandomness_2011} for the standard RSVD. Also, the term carrying the deviation from the optimum error depends on the target rank $k$, while it depends on $\sqrt{k}$ in~\cite{HalkoMartinssonEtAl_FindingStructureRandomness_2011}. This suggests that the generalization has been done, to some extents, to the detriment of the accuracy and that the bounds may overestimate the actual error.

    In this work, we propose a \textit{general stochastic error analysis} of the low-rank approximation error~\eqref{eq:metric} which holds for any matrix of the form $\ZBold = \A (\A\tra\A)^q \, \LBold \G$ where the columns of $\G$ are independently sampled from a standard multivariate Gaussian distribution.
    The obtained bounds highlight key properties that the matrix $\LBold$ should satisfy to achieve optimal low-rank approximation errors.
    Then, our general error analysis incorporates other existing stochastic error analyses, such as those in~\cite{HalkoMartinssonEtAl_FindingStructureRandomness_2011,Gu_SubspaceIterationRandomization_2015, BoulleTownsend_GeneralizationRandomizedSingular_2022}, for particular choices of $\LBold$ and $q$, thereby \textit{unifying a large variety of error analyses}.
    In particular, we show that choosing $\LBold = \I_n$ and $q=0$ yields exactly the same expectation bound as the one provided in ~\cite{HalkoMartinssonEtAl_FindingStructureRandomness_2011}. In this regard, the proposed generalization does not compromise the accuracy of the expectation bound for the particular choice of the identity.
    Additionally, our error analysis turns out to significantly enhance the accuracy of the bounds from~\cite{BoulleTownsend_GeneralizationRandomizedSingular_2022}. 
    Given the general form of $\ZBold$ considered, our error analysis can be applicable to a wide range of existing and randomized algorithms.
    
    The outline of this paper is the following.  We first introduce key background material that will be helpful throughout the manuscript in Section~\ref{sec:preliminary}. Section~\ref{sec:main} presents our main contribution on the general error analysis. In Section~\ref{sec:application}, we specialize our analysis to four practical cases: the RSVD $(\LBold = \I_n$ and $q=0$), the power iteration scheme ($\LBold = \I_n$ and $q \geq 0$), the generalized RSVD~\cite{BoulleTownsend_GeneralizationRandomizedSingular_2022} ($q = 0$) and the case where the covariance matrix is constructed out of a priori information on $\A$. To illustrate the versatility of Theorem~\ref{th:main}, we propose for the first three cases a comparison with~\cite{HalkoMartinssonEtAl_FindingStructureRandomness_2011,Gu_SubspaceIterationRandomization_2015,BoulleTownsend_GeneralizationRandomizedSingular_2022} respectively. In Section~\ref{sec:num}, we propose numerical illustrations on a data assimilation problem~\cite{Daley_AtmosphericDataAnalysis_1999,LahozKhattatovEtAl_DataAssimilation_2010} and show that exploiting the particular structure of the problem to define the covariance matrix $\CBold$ improves the performance of the RSVD. Finally, Section~\ref{Sec:conclusions} presents the conclusions together with the outlook.

\section{Preliminaries} \label{sec:preliminary}

    In this section, we recall well-known key results and definitions from the literature.

    \paragraph{Submultiplicativity}
        Let $\norm{\cdot}_2$ and $\norm{\cdot}_F$ denote the spectral and the Frobenius norm respectively. The strong submultiplicativity property~\cite[Relation (B.7)]{Higham_FunctionsMatricesTheory_2008} reads
        \begin{equation} \label{eq:submultiplicativity}
            \forall~\MBold \in \R^{n \times p}, \NBold \in \R^{p \times q}, \quad
            \norm{\MBold \NBold}_F \leq \norm{\MBold}_F \norm{\NBold}_2.
        \end{equation}

    \paragraph{Partial ordering on the set of symmetric matrices}
        Let $\MBold \in \Rmm, \NBold \in\Rmm$ be two symmetric matrices. The notation $\MBold \lleq \NBold$ means that $\NBold - \MBold$ is positive semi-definite. This relation defines a partial ordering on the set of symmetric matrices~\cite[Section 7.7]{HornJohnson_MatrixAnalysis_2012}. An important property~\cite[Theorem 7.7.2]{HornJohnson_MatrixAnalysis_2012} is that the partial ordering is preserved under the conjugation rule, i.e.,
        \begin{equation} \label{eq:conjugation_rule}
            \MBold \lleq \NBold ~\implies~
            \QBold\tra \MBold \QBold \lleq \QBold\tra \NBold \QBold,
            \quad \forall~\QBold \in \Rmn.
        \end{equation}
        We note that as a consequence of~\cite[Corollary 7.7.4(c)]{HornJohnson_MatrixAnalysis_2012}, the trace is monotonic with respect to the partial ordering, i.e.,
        \begin{equation} \label{eq:monotonic}
            \MBold \lleq \NBold ~\implies~ \tr(\MBold) \leq \tr(\NBold).
        \end{equation}
        
    \paragraph{Projection matrices}
        Suppose that $\MBold \in \Rmn$ has full column rank with column range space denoted by $\range(\MBold)$. We denote by  $\MBold\pinv$ the left multiplicative inverse of $\MBold$, i.e., the Moore-Penrose inverse of $\MBold$, see, e.g.,~\cite{HornJohnson_MatrixAnalysis_2012}. The orthogonal projection on $\range(\MBold)$ is then given by $\Proj(\MBold) = \MBold \MBold\pinv$, in particular, one has $\range(\Proj(\MBold)) = \range(\MBold)$.
    
    \paragraph{Sherman-Morrison formula}
        Let $\MBold \in \Rmn$ and $\NBold \in\Rnm$ such that $\I_n + \NBold \MBold$ is non-singular. Then, $\I_m + \MBold \NBold$ is also non-singular and one has~\cite[Section 2.1.4]{GolubVanLoan_MatrixComputations_1996},
        \begin{equation} \label{eq:sherman}
            (\I_m + \MBold \NBold)\inv = \I_m - \MBold(\I_n + \NBold \MBold)\inv \NBold.
        \end{equation}
    
    \paragraph{Principal angles between subspaces}
        Let $\sub{M}, \sub{N} \subset \Rm$ be two $k$-dimensional subspaces, and $\MBold, \NBold \in \Rmk$ be two associated matrices with orthogonal columns satisfying $\range(\MBold) = \sub{M}$ and $\range(\NBold) = \sub{N}$ respectively. If $\sigma_i(\cdot)$ denotes the $i$-th largest singular value of a given matrix, then the principal angles between $\sub{M}$ and $\sub{N}$, denoted by $\theta_1, \, \dots, \, \theta_k$ are defined as follow~\cite{KnyazevArgentati_PrincipalAnglesSubspaces_2002},
        \begin{equation*}
            \theta_i = \arccos(\sigma_i(\MBold\tra \NBold)), \quad 1 \leq i \leq k.
        \end{equation*}
        Further details on this notion can be found in~\cite{PaigeWei_HistoryGeneralityCS_1994}.
        
    \paragraph{Tangent matrix of the principal angles}
        Let $[\, \MBold, \, \overline{\MBold}\, ] \in \Rnn$ be an orthogonal matrix with $\MBold \in \Rnk$, and $\NBold \in \Rnk$ a full column rank matrix. Then the following singular value decomposition~\cite[Theorem 3.1]{ZhuKnyazev_AnglesSubspacesTheir_2013} yields
        \begin{equation} \label{eq:tangent}
            \tanmat\left(\MBold, \, \NBold\right) \equiv
            \overline{\MBold}\tra \NBold (\MBold\tra \NBold)\pinv =
            \U \diag\left(\tan(\theta_1), \, \dots, \, \tan(\theta_k)\right) \V\tra,
        \end{equation}
        where $\U \in \R^{(n-k) \times k}$ has orthonormal columns, $\V \in \Rkk$ is orthogonal and scalars $\theta_1, \, \dots, \, \theta_k$ are the principal angles between $\range(\MBold)$ and $\range(\NBold)$.
        
\section{General randomized low-rank approximation error bounds} \label{sec:main}    
    In this section, we first present a stochastic error analysis of \eqref{eq:metric} for the case where $\ZBold$ is an arbitrary $m \times \ell$ matrix (see Theorem~\ref{th:main:general}). More precisely, we consider the general case $\ZBold = \D \G \in \Rml$ where $\D \in \Rmr$ is an arbitrary matrix and $\G \in \Rrl$ is a matrix whose columns are independently sampled from a standard $r$-variate Gaussian distribution such that $r \geq \ell$ (to ensure $\G$ is full column rank). By construction, the columns of $\ZBold$ are drawn from a Gaussian distribution with zero mean and covariance matrix $\K = \D\D\tra \in \Rmm$. This result then serves as the foundation for our main stochastic error analysis, presented in Theorem~\ref{th:main}, which considers matrices of the form $\ZBold = \A(\A\tra\A)^q \, \LBold\G$. Algorithm \ref{alg:grsvd} presents the pseudo-code of a general randomized procedure that can typically be theoretically analyzed by Theorem~\ref{th:main}.
    \begin{algorithm}
        \caption{A general randomized singular value decomposition framework.}
        \label{alg:grsvd}
        \begin{algorithmic}[1]
                \Require{Matrices $\A \in \Rmn$ and $\LBold \in \Rnr$, number of samples $\ell \leq \min\set{\rank(\A), \rank(\LBold)}$, target rank $k \leq \ell$, and number of power iteration steps $q\geq 0$.}
                
                \State{Compute $\ZBold = \A \LBold \G \in \Rnl$  with $\G \in \R^{r \times \ell}$ whose columns are standard Gaussian vectors}
                \State{Perform the thin QR factorization $\ZBold = \QBold \RBold$}

                \For{$i = 0, \dots, q$}
                    \State{Compute $\ZBold = \A\tra\A \QBold $}
                    \State{Perform the thin QR factorization $\ZBold = \QBold \RBold$}
                \EndFor
            
                \State{Compute $\Y = [\A\QBold]\tra$}
                \State{Perform the order $k$ truncated SVD of $\Y$, i.e. $\Y = \U_k \widehat{\SigmaBold}_k \widehat{\V}_k\tra$}
                \State{Compute $\widehat{\U}_k = \QBold \U_k$}
                
                \Return{Matrices $\widehat{\U}_k \in \Rnk$, $\widehat{\SigmaBold}_k \in \Rkk$ and $\widehat{\V}_k \in \Rnk$ such that $\A \approx \widehat{\U}_k \widehat{\SigmaBold}_k \widehat{\V}_k\tra$.}
            \smallskip
        \end{algorithmic}
    \end{algorithm}
    
    Let us now state the general stochastic error analysis. For readability, 
    the proof of Theorem \ref{th:main:general} is given in Appendix~\ref{app:proofs}.

    \begin{theorem} \label{th:main:general}
        Let $\A \in \Rmn$ be an arbitrary matrix, and $\U_k$, $\SigmaBold_k$ and $\overline{\SigmaBold}_k$ the matrices related to the singular value decomposition of $\A$ and defined in \eqref{eq:partitioning_svd}. Let $\D \in \Rmr$ be a matrix satisfying $\rank(\D\tra \U_k) = k$, and let us define $\ZBold = \D \G \in \Rml$ with $\G \in \Rrl$ a matrix whose columns are independently sampled from an $r$-variate standard Gaussian distribution such that $\ell \leq \min\set{\rank(\A), \, \rank(\D)} \leq r$.
        
        \vspace{0.25cm}

        \noindent Then, for any integer $\ell \geq k + 2$, 
        \begin{equation*} 
        \label{th1:eq:1}\expect{\norm{[\I_n - \Proj(\ZBold)] \A}_F} \leq
            \left(1 + \tau_k^2 + \frac{\rho_k^2}{\ell - k - 1} \right)\half \norm{\overline{\SigmaBold}_k}_F.
        \end{equation*}

        \noindent Moreover, if $\ell \geq k + 4$, then for all $u, t \geq 1$,
        \begin{equation*}
        \label{th1:eq:2}    \norm{[\I_n - \Proj(\ZBold)] \A}_F \leq
            \left(1 + \tau_k + \sqrt{3}ut \cdot \frac{\rho_k}{\sqrt{\ell - k + 1}} \right) \norm{\overline{\SigmaBold}_k}_F
        \end{equation*}
        holds with probability at least $1 - e^{-u^2 / 2} - t^{-(\ell - k)}$ with 
        \begin{align*}
            & \tau_k = \frac{\norm{\tanmat(\U_k, \, \D\D\tra \U_k) \SigmaBold_k}_F}{\norm{\overline{\SigmaBold}_k}_F} \,
             \quad \txtand \quad
             \rho_k = \frac{  \norm{[\I_r - \Proj(\D\tra \U_k)] \D\tra}_F \,
            \norm{\SigmaBold_k (\U_k\D\D\tra\U_k)\halfinv}_F}{\norm{\overline{\SigmaBold}_k}_F} \,
          .
        \end{align*}
    \end{theorem}

    As presented in introduction, randomized algorithms available in the literature mostly rely on the particular choice $\D = \A(\A\tra\A)^q \, \LBold$, for $q \geq 0$ and $\LBold \in \Rnr$. In this particular case the coefficients $\tau_k$ and $\rho_k$ in Theorem~\ref{th:main:general} can be simplified, as stated in the next result. In Section~\ref{sec:application}, we show that this result unifies a wide variety of existing stochastic low-rank approximation error analyses. 
    
    \begin{theorem} \label{th:main}
        Let $\A \in \Rmn$ be an arbitrary matrix, and $\V_k$, $\overline{\V}_k$, $\SigmaBold_k$ and $\overline{\SigmaBold}_k$ the matrices related to the singular value decomposition of $\A$ and defined in \eqref{eq:partitioning_svd}. Let $\LBold \in \Rnr$ be a matrix satisfying $\rank(\LBold\tra \V_k) = k$, and let us define $\ZBold = \A(\A\tra\A)^q \, \LBold \G \in \Rml$ with $q\geq 0$ and $\G \in \Rrl$ a matrix whose columns are independently sampled from an $r$-variate standard Gaussian distribution such that $\ell \leq \min\set{\rank(\A), \, \rank(\LBold)} \leq r$. Then the expression of $\tau_k$ and $\rho_k$ in Theorem \ref{th:main:general} simplifies to,
        \begin{align*}
            \tau_k & =
            \frac{1}{\norm{\overline{\SigmaBold}_k}_F} \,
            \norm{(\overline{\SigmaBold}_k \overline{\SigmaBold}_k\tra)^q \overline{\SigmaBold}_k \tanmat(\V_k, \, \LBold\LBold\tra \V_k) \SigmaBold_k^{-2q}}_F, \\
            \txtand \quad \rho_k & =
            \frac{1}{\norm{\overline{\SigmaBold}_k}_F} \,
            \norm{[\I_r - \Proj(\LBold\tra \V_k)] \LBold\tra \overline{\V}_k \overline{\SigmaBold}_k\tra \, (\overline{\SigmaBold}_k \, \overline{\SigmaBold}_k\tra)^q}_F \,
            \norm{\SigmaBold_k^{-2q} (\V_k\tra \LBold\LBold\tra \V_k)\halfinv}_F.
        \end{align*}
        In this case, $\tau_k$ and $\rho_k $ can be bounded as follows:
    \begin{equation*}
            \tau_k \leq \norm{\tanmat(\V_k, \, \LBold\LBold\tra \V_k)}_2 \left(\dfrac{\sigma_{k+1}}{\sigma_k}\right)^{2q}
            \quad \txtand \quad \rho_k \leq \kappa_2(\LBold) \, \sqrt{k} \left(\dfrac{\sigma_{k+1}}{\sigma_k}\right)^{2q},
        \end{equation*}
        where $\kappa_2(\LBold)$ denotes the 2-norm condition number of the matrix $\LBold$.
    \end{theorem}
    \begin{proof}.~
        First, note that from $\A(\A\tra\A)^q = \U \SigmaBold (\SigmaBold\tra \SigmaBold)^q \V\tra$, we get
        \begin{equation} \label{eq:d_tra_u}
            \D\tra \U_k =
            \LBold\tra \V_k \SigmaBold_k^{2q+1} \quad \txtand \quad
            \D\tra \overline{\U}_k =
            \LBold\tra \overline{\V}_k (\overline{\SigmaBold}_k\tra \overline{\SigmaBold}_k)^q \overline{\SigmaBold}_k\tra.
        \end{equation}
        Since $k \leq \rank(\A)$, then $\SigmaBold_k$ is full rank and $\rank(\D\tra \U_k) = \rank(\LBold\tra \V_k) = k$. Then, one has $\rank(\D) = \min\set{\rank(\A), \, \rank(\LBold)}$ by definition of $\D$ which implies that $\min\set{\rank(\A), \, \rank(\D)} = \min\set{\rank(\A), \, \rank(\LBold)}$. Consequently, one can apply Theorem \ref{th:main:general}.

        Let us now prove the simplifications in the coefficients. For $\tau_k$, using \eqref{eq:d_tra_u} yields the following,
        \begin{align*}
            \tanmat(\U_k, \, \D\D\tra \U_k) \SigmaBold_k & =
            \overline{\U}_k\tra \D\D\tra \U_k (\U_k\tra \D\D\tra \U_k)\inv \SigmaBold_k \\ & =
            \overline{\SigmaBold}_k (\overline{\SigmaBold}_k\tra \overline{\SigmaBold}_k)^q \overline{\V}_k\tra \LBold\LBold\tra \V_k (\V_k\tra \LBold\LBold\tra \V_k)\inv \SigmaBold_k^{-2q}.
        \end{align*}
        Identifying that $\tanmat(\V_k, \, \LBold\LBold\tra \V_k) = \overline{\V}_k\tra \LBold\LBold\tra \V_k (\V_k\tra \LBold\LBold\tra \V_k)\inv$, one readily gets the result for $\tau_k$. For $\rho_k$, one has from $\range(\D\tra \U_k) = \range(\LBold\tra \V_k)$ that $\Proj(\D\tra \U_k) = \Proj(\LBold\tra \V_k)$, and then using~\eqref{eq:d_tra_u} we obtain,
        \begin{equation*}
            [\I_r - \Proj(\LBold\tra \V_k)] \D\tra \overline{\U}_k\tra =
            [\I_r - \Proj(\LBold\tra \V_k)] \LBold \overline{\V}_k
            (\overline{\SigmaBold}_k\tra \overline{\SigmaBold}_k)^q \overline{\SigmaBold}_k\tra \overline{\U}_k\tra.
        \end{equation*}
        Finally, we use the following relation,
        \begin{equation*}
            \SigmaBold_k (\U_k\tra \D\D\tra \U_k)\inv \SigmaBold_k =
            \SigmaBold_k (\SigmaBold_k^{2q+1} \V_k\tra \LBold\LBold\tra \V_k \SigmaBold_k^{2q+1})\inv \SigmaBold_k =
            \SigmaBold_k^{-2q} (\V_k\tra \LBold\LBold\tra \V_k)\inv \SigmaBold_k^{-2q},
        \end{equation*}
        to get the simplified expressions for $\tau_k$ and $\rho_k$.

         Now,  for the bounds on $\tau_k$ and $\rho_k$, we proceed as follows. First, for $\tau_k$, we use the norm submultiplicativity property to get  
        \begin{align*}
            \norm{(\overline{\SigmaBold}_k \overline{\SigmaBold}_k\tra)^q \overline{\SigmaBold}_k \tanmat(\V_k, \, \LBold\LBold\tra \V_k) \SigmaBold_k^{-2q}}_F & \leq \norm{(\overline{\SigmaBold}_k\tra \overline{\SigmaBold}_k)^q}_2 \norm{\overline{\SigmaBold}_k}_F \norm{\tanmat(\V_k, \, \LBold\LBold\tra \V_k)}_2 \norm{\SigmaBold_k^{-2q}}_2 \\
            & = \norm{\overline{\SigmaBold}_k}_F\, \sigma_{k+1}^{2q} \, \norm{\tanmat(\V_k, \, \LBold\LBold\tra \V_k)}_2 \, \sigma_{k}^{-2q},
        \end{align*}
        from which we readily deduce the bound for $\tau_k$. 
        
        Then, regarding  the bound on $\rho_k$, as for $\tau_k$, we use the  norm submultiplicativity to get,
        \begin{align*}
            \norm{[\I_r - \Proj(\LBold\tra \V_k)] \LBold\tra \overline{\V}_k \overline{\SigmaBold}_k\tra \, (\overline{\SigmaBold}_k \, \overline{\SigmaBold}_k\tra)^q}_F & \leq 
            \norm{[\I_r - \Proj(\LBold\tra \V_k)]}_2
            \norm{\LBold\tra}_2
            \norm{\overline{\V}_k}_2
            \norm{\overline{\SigmaBold}_k\tra}_F \, \norm{(\overline{\SigmaBold}_k \, \overline{\SigmaBold}_k\tra)^q}_2 \\
            & = \lambda_1(\LBold) \, \norm{\overline{\SigmaBold}_k}_F \, \sigma_{k+1}^{2q},
        \end{align*}
        where $\lambda_1(\cdot)$ denotes the largest eigenvalue. Then, since $\V_k$ is orthogonal, one has $(\V_k\tra \LBold\LBold\tra \V_k) \succcurlyeq \lambda_n(\LBold\LBold\tra) \I_n = \lambda_n(\LBold)^2 \I_n$, where $\lambda_n(\cdot)$ denotes the smallest eigenvalue. Using the conjugation rule \eqref{eq:conjugation_rule} yields $\SigmaBold_k^{-2q} (\V_k\tra \LBold\LBold\tra \V_k)\inv \SigmaBold_k^{-2q} \lleq \frac{1}{\lambda_n(\LBold)^2} \SigmaBold_k^{-4q}$ and from the trace monotonicity \eqref{eq:monotonic} one has, 
        \begin{equation*}
            \norm{\SigmaBold_k^{-2q} (\V_k\tra \LBold\LBold\tra \V_k)\halfinv}_F^2=
            \tr\left(\SigmaBold_k^{-2q} (\V_k\tra \LBold\LBold\tra \V_k)\inv \SigmaBold_k^{-2q}\right) \leq
            \dfrac{\tr(\SigmaBold_k^{-4q})}{\lambda_n(\LBold)^2} =
            \dfrac{\norm{\SigmaBold_k^{-2q}}_F}{\lambda_n(\LBold)^2} \leq
            k\dfrac{\sigma_k^{-2q}}{\lambda_n(\LBold)^2}.
        \end{equation*}
        Noting that $\kappa_2(\LBold) = \lambda_1(\LBold) / \lambda_n(\LBold)$, altogether, we obtain the desired inequality for $\rho_k$.
    \end{proof}

    Theorem \ref{th:main} predicts that the deviation from the optimal error depends on two coefficients, namely $\tau_k$ and $\rho_k$.
   The bounds in $\tau_k$ and $\rho_k$, as given in Theorem~\ref{th:main}, also illustrate how the error behaves with respect to the power iteration parameter~$q$. They indicate that, regardless of~$\LBold$, and as long as $\sigma_{k+1} < \sigma_k$, the constants~$\tau_k$ and~$\rho_k$ get smaller as~$q$ increases.
    
    Theorem~\ref{th:main} also suggests that the ideal choice of $\LBold$ must yield $\tau_k = \rho_k = 0$.
    In fact, having $\tau_k = 0$ is equivalent to $\tanmat(\V_k, \, \LBold\LBold\tra \V_k) = \bm{0}$, that is $\range(\LBold\LBold\tra \V_k) = \range(\V_k)$. This condition is fulfilled whenever $\range(\V_k)$ is invariant under the action of $\LBold\LBold\tra$, or equivalently whenever $\range(\V_k)$ is an eigenspace of $\LBold\LBold\tra$. For $\rho_k$, let us illustrate its behavior on the following particular symmetric case,
    \begin{equation*}
        \LBold = \V_k\V_k\tra + \beta \overline{\V}_k \overline{\V}_k\tra, \quad \beta > 0.
    \end{equation*}
    Here, $\LBold\tra\V_k = \V_k$ and $\LBold\tra\overline{\V}_k = \beta \overline{\V}_k$, from which we trivially deduce 
    \begin{equation*}
        \rho_k = \beta
        \frac{\norm{\overline{\SigmaBold}_k\tra \, (\overline{\SigmaBold}_k \, \overline{\SigmaBold}_k\tra)^q}_F}{\norm{\overline{\SigmaBold}_k}_F} \,
        \norm{\SigmaBold_k^{-2q}}_F.
    \end{equation*}
    Consequently, $\rho_k$ is small if $\beta$ is small too, that is if the dominant eigenspace of $\LBold$ is carried by $\V_k$, and the remaining part of the spectrum is negligible. In this sense, $\rho_k$ can be interpreted as a separation coefficient between $\range(\V_k)$ and its orthogonal. Optimally, if $\beta = 0$, then $\range(\LBold) = \range(\V_k)$ and $\rho_k = 0$.

       Consequently, an optimal choice for $\LBold$ would be of the form,
        \begin{equation} \label{eq:opt_L}
            \LBold = \V_k \X, \quad \text{with} \quad \X \in \mathbb{R}^{k \times r}.
        \end{equation}
        In practice, since $\V_k$ is precisely the quantity that needs to be determined, which makes~\eqref{eq:opt_L} out of reach. Nevertheless, by observing that $\A\tra = \A_k\tra + \overline{\A}_k\tra$ and the fact that $\A_k\tra = \V_k \SigmaBold_k \U_k\tra$ is of the same form as \eqref{eq:opt_L}, then the bounds predict that the RSVD will be particularly efficient if $\A\tra \approx \A_k\tra$, that is whenever $\overline{\A}_k\tra$ is negligible compared to $\A_k\tra$. We recover here a property that is well documented in the randomized linear algebra literature.

\section{Particular choices for the covariance matrix} \label{sec:application}

    In this section, we particularize Theorem~\ref{th:main} for certain choices of $\LBold$ and $q$. As we will show, several prior bounds, that have been independently derived, can be recovered from Theorem~\ref{th:main}. In Section \ref{sec:application:rsvd}, we propose a stochastic analysis of the standard RSVD that can be put into perspective with the  reference analysis in~\cite{HalkoMartinssonEtAl_FindingStructureRandomness_2011}. Secondly, Section \ref{sec:application:power_iteration} provides bounds for the randomized power iteration that are compared to the ones in~\cite{Gu_SubspaceIterationRandomization_2015}. Then, in Section \ref{sec:application:grsvd} we propose an analysis of the Generalized RSVD introduced in~\cite{BoulleTownsend_GeneralizationRandomizedSingular_2022} which highlights the gain in accuracy brought by Theorem \ref{th:main}. Finally, we discuss in Section \ref{sec:application:available_approx} a relevant choice for the covariance matrix $\LBold$ when a priori information on $\A$ is available.

    \subsection{The RSVD ($\LBold = \I_n$ and $q=0$)} \label{sec:application:rsvd}
        The RSVD corresponds to the special case where $\LBold = \I_n$ and $q=0$. In this case, the bounds of Theorem~\ref{th:main} simplify, as shown in Corollary \ref{cor:hmt}.
        \begin{corollary} \label{cor:hmt}
            Let $\A \in \Rmn$ be an arbitrary matrix, and $\overline{\SigmaBold}_k$ the matrix related to the singular value decomposition of $\A$ and defined in \eqref{eq:partitioning_svd}. Let us define $\ZBold = \A \G \in \Rml$ with $\G \in \Rnl$ a matrix whose columns are independently sampled from an $n$-variate standard Gaussian distribution such that $\ell \leq \rank(\A)$.
            
            \vspace{0.25cm}

            \noindent Then, for any integer $\ell \geq k + 2$, one has
            \begin{equation}
                \expect{\norm{[\I_n - \Proj(\ZBold)] \A}_F} \leq
                \left(1 + \sqrt{\frac{k}{\ell - k - 1}} \right) \norm{\overline{\SigmaBold}_k}_F.
            \end{equation}
            
            \noindent Moreover, if $\ell \geq k + 4$, then for all $u, t \geq 1$,
            \begin{equation}
                \norm{[\I_n - \Proj(\ZBold)] \A}_F \leq
                \left(1 + \sqrt{3} ut \cdot \sqrt{\frac{k}{\ell - k + 1}} \right) \norm{\overline{\SigmaBold}_k}_F,
            \end{equation}
            holds with probability at least $1 - e^{-u^2 / 2} - t^{-(\ell - k)}$.    
        \end{corollary}
        \begin{proof}.~
            This corollary is a consequence of Theorem \ref{th:main} when $q=0$ and $\LBold = \I_n$, where firstly we observe that $\tanmat(\V_k, \, \LBold\LBold\tra\V_k) = \tanmat(\V_k, \, \V_k) = \bf{0}$, and secondly that $\kappa_2(\LBold) = 1$, from which we immediately gets $\tau_k(\I_n) = 0$ and $\rho_k(\I_n) \leq \sqrt{k}$.
        \end{proof}
        Corollary \ref{cor:hmt} provides the same expectation bound as in~\cite[Theorem 10.5]{HalkoMartinssonEtAl_FindingStructureRandomness_2011}, and offers comparable bounds in probability. For easier comparison of the bounds in probability, we recall~\cite[Theorem 10.7] {HalkoMartinssonEtAl_FindingStructureRandomness_2011} which states that under the same assumptions as Corollary~\ref{cor:hmt}, 
        \begin{equation}
            \norm{[\I_n - \Proj(\ZBold)] \A}_F \leq \left(1 + \sqrt{3} t \cdot \sqrt{\frac{k}{\ell - k + 1}} \right) \norm{\overline{\SigmaBold}_k}_F + ut \frac{e \sqrt{\ell}}{\sqrt{\ell-k+1}}\sigma_{k+1}
        \end{equation}
        holds with probability at least $1 - e^{-u^2 / 2} - 2t^{-(\ell - k)}$. Our bound is almost identical to the first term in $\norm{\overline{\SigmaBold}_k}_F$ except for the extra factor $u$. However, since our bound does not have the second term in $\sigma_{k+1}$, then for reasonably small values of $u$ (\eg 2 or 3) one can expect our probability bound to have comparable accuracy.

    \subsection{The randomized power iteration ($\LBold = \I_n$ and $q \geq 0$)} \label{sec:application:power_iteration}
        The power iteration scheme considers $\ZBold = \A (\A\tra\A)^q \G$, that is $\LBold = \I_n$ and $q \geq 0$. Exploiting this structure in Theorem~\ref{th:main} yields the following corollary.

        \begin{corollary} \label{cor:gu}
            Let $\A \in \Rmn$ be an arbitrary matrix, and $\overline{\SigmaBold}_k$ the matrix related to the singular value decomposition of $\A$ and defined in \eqref{eq:partitioning_svd}. Let us define $\ZBold = \A (\A\tra \A)^q \G$ with $q \geq 0$ and $\G \in \Rnl$ a matrix whose columns are independently sampled from an $n$-variate standard Gaussian distribution such that $\ell \leq \rank(\A)$.
            
            \vspace{0.25cm}

            \noindent Then, for any integer $\ell \geq k + 2$, one has
            \begin{equation}
                \expect{\norm{[\I_n - \Proj(\ZBold)] \A}_F} \leq
                \left(1 + \left(\frac{\sigma_{k+1}}{\sigma_k}\right)^{2q} \sqrt{\frac{k}{\ell - k - 1}} \right) \norm{\overline{\SigmaBold}_k}_F.
            \end{equation}
            
            \noindent Moreover, if $\ell \geq k + 4$, then for all $u, t \geq 1$,
            \begin{equation}
                \norm{[\I_n - \Proj(\ZBold)] \A}_F \leq
                \left(1 + \sqrt{3} ut \cdot \left(\frac{\sigma_{k+1}}{\sigma_k}\right)^{2q} \sqrt{\frac{k}{\ell - k + 1}} \right) \norm{\overline{\SigmaBold}_k}_F,
            \end{equation}
            holds with probability at least $1 - e^{-u^2 / 2} - t^{-(\ell - k)}$.    
        \end{corollary}
        \begin{proof}.~
            This corollary is also a consequence of Theorem~\ref{th:main} when $\LBold = \I_n$. We have already seen that $\tau_k(\I_n) = 0$, and observing that $\kappa_2(\I_n) = 1$ we get $\rho_k(\I_n) \leq \sqrt{k} \, (\sigma_{k+1} / \sigma_k)^{2q}$.            
        \end{proof}
    
        For comparison, using the property that $\sqrt{1 + x} \leq 1 + \sqrt{x}$ for $x \geq 0$, Theorem 5.7 in~\cite{Gu_SubspaceIterationRandomization_2015} provides, under the same assumptions as Corollary \ref{cor:gu}, the following upper bound:
        \begin{equation*}
            \expect{\norm{[\I_n - \Proj(\ZBold)] \A}_F} \leq
            \left(1 +
            \sqrt{k}\left(\frac{\sigma_{k+1}}{\sigma_k}\right)^{2q} \frac{4e \, \sqrt{\ell}}{\ell - k + 1} \left[\sqrt{n-k} + \sqrt{\ell} + 7\right] \frac{\sigma_{k+1}}{\norm{\overline{\SigmaBold}_k}_F} \right)
            \norm{\overline{\SigmaBold}_k}_F.
        \end{equation*}
        This bound is similar to our bound except for the multiplying factor
        \begin{equation*}
            \mathsf{c} = 4e \, \sqrt{\ell} \, \frac{\sqrt{\ell - k - 1}}{\ell - k + 1} \left[\sqrt{n-k} + \sqrt{\ell} + 7\right] \frac{\sigma_{k+1}}{\norm{\overline{\SigmaBold}_k}_F}.
        \end{equation*}
        Since $\norm{\overline{\SigmaBold}_k}_F \leq \sigma_{k+1} \, \sqrt{n-k}$ from classical norm equivalence, and both $\ell \geq \ell - k - 1$ and $\ell - k \geq 2$ by assumption, it yields
        \begin{equation*}
            \mathsf{c} \geq \frac{4e}{3} \, (\ell - k - 1).
        \end{equation*}
        This suggests that Corollary~\ref{cor:gu} provides tighter bounds than the ones proposed in~\cite{Gu_SubspaceIterationRandomization_2015}, especially when the oversampling $\ell - k$ is large. Nevertheless, the refinement might not be significant and only illustrates that bounds from~\cite{Gu_SubspaceIterationRandomization_2015} are slightly over-pessimistic.
    
    \subsection{The generalized RSVD ($q=0$)} \label{sec:application:grsvd}
        Let us now consider $\ZBold = \A \LBold \G$, that is $q=0$. This configuration has been studied in~\cite{BoulleTownsend_GeneralizationRandomizedSingular_2022}, to improve the accuracy of low-rank approximations of differential operators. Our corresponding analysis obtained from Theorem \ref{th:main} is proposed in Corollary \ref{cor:bt}.

        \begin{corollary} \label{cor:bt}
            Let $\A \in \Rmn$ be an arbitrary matrix, and $\V_k$, $\overline{\V}_k$, $\SigmaBold_k$ and $\overline{\SigmaBold}_k$ the matrices related to the singular value decomposition of $\A$ and defined in \eqref{eq:partitioning_svd}. Let $\LBold \in \Rnr$ be a matrix satisfying $\rank(\LBold\tra \V_k) = k$ and let us define $\ZBold = \A \LBold \G \in \Rml$ with $\G \in \Rrl$ a matrix whose columns are independently sampled from an $r$-variate standard Gaussian distribution such that $\ell \leq \min\set{\rank(\A), \, \rank(\LBold)} \leq r$. Let us further denote $\CBold = \LBold\LBold\tra$.
            
            \vspace{0.25cm}

            \noindent Then, for any integer $\ell \geq k + 2$,
            \begin{equation*}
                \expect{\norm{[\I_n - \Proj(\ZBold)] \A}_F} \leq
                \left(1 + \sqrt{\ell - k} \, \sqrt{\frac{k}{\ell-k-1} \, \frac{\beta_k}{\gamma_k}} \right) \norm{\overline{\SigmaBold}_k}_F.
            \end{equation*}
            
            \noindent Moreover, if $\ell \geq k + 4$, then for all $u, t \geq 1$,
            \begin{equation*}
                \norm{[\I_n - \Proj(\ZBold)] \A}_F \leq
                \left(1 + \sqrt{3} \, \left(\sqrt{\ell - k + 1} + 1 \right) \, \sqrt{\frac{k}{\ell-k+1} \, \frac{\beta_k}{\gamma_k}} \cdot ut \right) \norm{\overline{\SigmaBold}_k}_F,
            \end{equation*}
            holds with probability at least $1 - e^{-u^2 / 2} - t^{-(\ell - k)}$.

            \vspace{0.25cm}

            \noindent Here,
            \begin{equation*}
                \beta_k = \frac{\tr\left(\overline{\SigmaBold}_k\tra \overline{\SigmaBold}_k \overline{\V}_k\tra \CBold \overline{\V}_k \right)}{\lambda_1(\CBold) \, \norm{\overline{\SigmaBold}_k}_F^2} \quad \txtand \quad
                \gamma_k = \frac{k}{\lambda_1(\CBold) \, \tr\left((\V_k\tra \CBold \V_k)\inv \right)}.
            \end{equation*}
        \end{corollary}
        \begin{proof}.~
            The assumptions allows to apply Theorem \ref{th:main} with $q=0$, so it remains to deal with the coefficients $\tau_k$ and $\rho_k$. Using the submultiplicativity of the Frobenius norm along with~\eqref{eq:tangent}, we have
            \begin{align*}
                \norm{\overline{\SigmaBold}_k \tanmat(\V_k, \, \CBold \V_k)}_F^2
                & = \norm{\overline{\SigmaBold}_k \overline{\V}_k\tra \CBold \V_k (\V_k\tra \CBold \V_k)\inv}_F^2 \\
                & \leq \norm{\overline{\SigmaBold}_k \, \overline{\V}_k\tra \LBold}_F^2 \norm{\LBold\tra \V_k (\V_k\tra \CBold \V_k)\inv}_F^2 \\
                & = \tr\left(\overline{\SigmaBold}_k\tra \overline{\SigmaBold}_k \overline{\V}_k\tra \CBold \overline{\V}_k \right) \cdot \tr\left((\V_k\tra \CBold \V_k)\inv \right) = k \, \norm{\overline{\SigmaBold}_k}_F^2 \, \frac{\beta_k}{\gamma_k}.
            \end{align*}
            Thus,
            \begin{equation*}
                \tau_k =
                \frac{\norm{\overline{\SigmaBold}_k \tanmat(\V_k, \, \CBold \V_k)}_F}{\norm{\overline{\SigmaBold}_k}_F} \leq \sqrt{k} \, \sqrt{\frac{\beta_k}{\gamma_k}}
            \end{equation*}
            Then, let us use the following fact,
            \begin{align*}
                \norm{[\I_m - \Proj(\LBold \V_k)] \LBold \overline{\V}_k \overline{\SigmaBold}_k\tra}^2_F & =
                \tr\left(\overline{\SigmaBold}_k\overline{\V}_k\tra \CBold \overline{\V}_k\overline{\SigmaBold}_k\tra - \overline{\SigmaBold}_k\overline{\V}_k\tra \CBold \V_k \CBold_k\inv \V_k\tra \CBold \overline{\V}_k\overline{\SigmaBold}_k\tra\right).
            \end{align*}
            Since $\overline{\SigmaBold}_k\overline{\V}_k\tra \CBold \V_k \CBold_k\inv \V_k\tra \CBold \overline{\V}_k\overline{\SigmaBold}_k\tra$ is symmetric positive semi-definite, one has,
            \begin{equation*}
                \norm{[\I_m - \Proj(\LBold \V_k)] \LBold \overline{\V}_k \overline{\SigmaBold}_k\tra}^2_F \leq
                \tr\left(\overline{\SigmaBold}_k\overline{\V}_k\tra \CBold \overline{\V}_k\overline{\SigmaBold}_k\tra \right)
                = \tr\left(\overline{\SigmaBold}_k\tra \overline{\SigmaBold}_k \overline{\V}_k\tra \CBold \overline{\V}_k \right).
            \end{equation*}
            Altogether, this yields
            \begin{equation*}
                \rho_k \leq \sqrt{k} \, \sqrt{\frac{\beta_k}{\gamma_k}}.
            \end{equation*}
            Plugging these inequalities in Theorem~\ref{th:main} and using the property that $\sqrt{1 + x} \leq 1 + \sqrt{x}$ for all $x \geq 0$, provides the desired result.
        \end{proof}

        The bounds in Corollary~\ref{cor:bt} are purposely formulated using the same coefficients $\beta_k$ and $\gamma_k$ as in~\cite{BoulleTownsend_GeneralizationRandomizedSingular_2022}. This formulation clearly leads to a loss of accuracy compared to the bounds obtained by directly setting $q=0$ in Theorem \ref{th:main}, but it significantly simplifies the comparison. As a result, the bounds in Corollary~\ref{cor:bt} are similar to those in~\cite{BoulleTownsend_GeneralizationRandomizedSingular_2022}, with notable difference. The probability and expectation bounds in~\cite{BoulleTownsend_GeneralizationRandomizedSingular_2022}, stated respectively in Theorem 2 and Proposition 6, predict that the deviation from the optimal error behaves as $O(k)$, where $k$ is the target rank. In contrast, Corollary~\ref{cor:bt} shows that this deviation behaves as  $O(\sqrt{k})$. This shows that the influence of the target rank $k$ is not as penalizing as the bounds in~\cite{BoulleTownsend_GeneralizationRandomizedSingular_2022} would imply.
        
    \subsection{$\LBold$ as a low-rank matrix update} \label{sec:application:available_approx}
        Let us now explore the potential benefits of using $\ZBold = \A (\A\tra \A)^q \, \LBold \G$ where $\LBold$ is constructed out of available approximate SVD. Let $\widehat{\V}_k \in \Rnk$ be an approximation of $\V_k$ satisfying $\widehat{\V}_k\tra \widehat{\V}_k = \I_k$, and  $\widehat{\SigmaBold}_k = \diag(\hat{\sigma}_1, \, \dots, \hat{\sigma}_k)$ an approximation of $\SigmaBold_k$. Following~\eqref{eq:opt_L}, a natural choice to consider for $\LBold$ can be given by 
        \begin{equation*}
            \LBold = \widehat{\V}_k \widehat{\SigmaBold}_k \widehat{\V}_k\tra.
        \end{equation*}
        In this case, $\LBold$ is a rank $k$ matrix such that $\range(\LBold \V_k) = \range(\widehat{\V}_k)$, which trivially yields 
        \begin{equation*}
            \tau_k =
            \frac{1}{\norm{\overline{\SigmaBold}_k}_F} \,
            \norm{(\overline{\SigmaBold}_k \overline{\SigmaBold}_k\tra)^q \overline{\SigmaBold}_k \tanmat(\V_k, \, \widehat{\V}) \SigmaBold_k^{-2q}}_F
            \quad \txtand \quad
            \rho_k = 0.
        \end{equation*}
        In this configuration, the low-rank approximation error only depends on how accurately $\widehat{\V}_k$ approximates $\V_k$. Furthermore, since $\ZBold = \A (\A\tra \A)^q \, \LBold \G$, one also has $\range(\ZBold) = \range(\A (\A\tra \A)^q \widehat{\V}_k)$ with probability one. This means that this choice would actually result in a deterministic algorithm, where the range of $\A$ is obtained by power iteration with $\widehat{\V}_k$ as an initial guess.

        However, since $\widehat{\V}_k$ is only an approximation of $\V_k$, then $\V_k$ has a non-zero intersection with both $\range(\widehat{\V}_k)$ and its orthogonal. Consequently, if $\widehat{\V}_k$ is not an approximation with satisfying accuracy, the orthogonal of $\range(\widehat{\V}_k)$ still contains relevant information regarding $\V_k$, and should be explored using randomization. A simple way to achieve this consists in considering instead,
        \begin{equation*}
            \LBold_\beta =
            \widehat{\V}_k \widehat{\SigmaBold}_k \widehat{\V}_k\tra +
            \beta \hat{\sigma}_k (\I_n - \widehat{\V}_k \widehat{\V}_k\tra), \quad \text{with} \quad \beta > 0.
        \end{equation*}
        The presence of $\hat{\sigma}_k$ allows to consider $\beta$ as an dimensionless parameter. This normalized parameter $\beta$ allows to balance the confidence in the approximation $\widehat{\V}_k$. If $\beta$ is close to zero, then $\range(\LBold_\beta \G)$ would be mostly aligned with $\range(\widehat{\V}_k)$. Conversely, if $\beta$ is close to 1, then the second term of $\LBold_\beta$ has the same order of magnitude than the first one, and therefore $\range(\LBold_\beta \G)$ would instead be almost uniformly distributed. In this case, using $\LBold_\beta$ becomes irrelevant as it will not bring any additional information.

\section{Numerical experiments} \label{sec:num}
    We validate our theoretical results by applying randomized algorithms within the framework of variational data assimilation (VarDA)~\cite{Kalnay_2002, Bannister2017}, which is formulated as a nonlinear, regularized, weighted least-squares problem. Solving this nonlinear problem involves addressing a \textit{sequence of symmetric and positive definite (SPD) linear  systems}~\cite{GrattonLawlessEtAl_ApproximateGaussNewton_2007, NocedalWright_NumericalOptimization_2006}. Each linear system is solved by using an iterative method such as preconditioned conjugate gradient method (PCG)~\cite{Saad_IterativeMethodsLinear_2019, HestenesStiefel_MethodsConjugateGradients_1952}. The potential of randomized algorithms has recently been explored in various VarDA formulations. For example, \cite{BousserezGuerretteEtAl_EnhancedParallelizationIncremental_2020} proposed using randomized algorithms as an alternative to PCG method. In another study, \cite{DauzickaiteLawlessEtAl_RandomisedPreconditioningForcing_2021, DauzickaiteLawlessEtAl_RandomisedPreconditioningStrong_2021} employed randomized algorithms to estimate the eigenspectrum of an SPD matrix, which was then used to construct a preconditioner for different formulations of weak-constraint VarDA~\cite{Tremolet_weakconstraint_2006}. These studies primarily relied on randomized algorithms based on the standard Gaussian distribution, i.e., $\LBold = \I_n$.

    In our experiments, we aim to demonstrate the potential of using non-standard Gaussian distributions to reduce the error of the approximated matrix within a fixed computational budget. This is particularly important in VarDA, where computational resources are often tightly constrained. To this end, in Section~\ref{sec:3dvar} we present the  three-dimensional VarDA formulation along with the resulting SPD system. Section~\ref{sec:generalized_SVD_algo} introduces the generalized randomized algorithm employed to approximate the eigenspectrum of a SPD matrix. Section~\ref{subsec:accuracy_bounds} illustrates the accuracy of the bounds in Theorem~\ref{th:main}, compares these bounds with those proposed in~\cite{HalkoMartinssonEtAl_FindingStructureRandomness_2011, BoulleTownsend_GeneralizationRandomizedSingular_2022} together with the empirical error obtained by performing the algorithm outlined in Section~\ref{sec:generalized_SVD_algo}. Then, in Section \ref{subsec:choices_for_L} we study the potential of different choices for $\LBold$ to result in better performance. 

    In this section, we focus on the stochastic error analysis of the approximated SPD matrix obtained using randomized algorithms based on a non-standard Gaussian distribution. This general approach is versatile and can be employed either within an iterative randomized algorithm or to extract approximate information for use in a preconditioner. Analyzing the performance of the iterative solvers with this approach is beyond the scope of the current paper and is left for future work.
    
    \subsection{Three dimensional VarDA} \label{sec:3dvar}
        For our numerical experiments, we use a three-dimensional variational data assimilation formulation. In this approach, we solve a weighted nonlinear least-squares problem for a fixed time:
         \begin{equation} \label{eq:nlsq}
            \min_{\mathbf{x}\in\Rn} f(\mathbf{x}) = \frac{1}{2} \norm{\mathbf{y} - \mathcal{H}(\mathbf{x})}_{\RBold\inv}^2 + \frac{1}{2} \norm{\mathbf{x} - \mathbf{x}_b}_{\B\inv}^2,
        \end{equation}
        where $\mathbf{x} \in \Rn$ is the state of a dynamical system, for instance temperature, $\mathbf{y} \in \Rm$ is a vector consisting of observations, and $\mathbf{x}_b \in \Rn$ is a priori information. The operator $\mathcal{H}: \Rn \rightarrow \Rm$ is the nonlinear observation operator mapping the state vector from the model space to the observation space. $\B \in \Rnn$ and $\RBold \in \Rmm$ represent the error covariance matrices of a priori information and observations respectively.
        
        A common approach for solving~\eqref{eq:nlsq} is the truncated Gauss-Newton method~\cite{GrattonLawlessEtAl_ApproximateGaussNewton_2007}. This method relies on the linearization of the nonlinear observation operator around the current iterate $\mathbf{x}_j$, which results in a weighted linear least-squares sub-problem whose solution can be obtained by solving a preconditioned linear system:
        \begin{equation} \label{eq:da_linear_system}
            (\I_n + \W \H_j\tra \RBold\inv \H_j \W ) \delta \mathbf{v}_j = - \W \mathbf{b}_j
        \end{equation}
        where $\W \in \Rnn$ is a symmetric square-root factorization of $\B$ used as a preconditioner, $\H_j \in \Rmn$ is the Jacobian matrix of the observation operator at the current iterate, $\delta \mathbf{v}_j = \W \delta \mathbf{x}_j$ is the increment (search direction) in the preconditioned space, and $\mathbf{b}_j$ is the gradient of the nonlinear cost function~(\ref{eq:nlsq}) calculated at the current iterate. The linear system~(\ref{eq:da_linear_system}) consisting of the symmetric positive definite matrix,
        \begin{equation} \label{eq:da_matrix}
            \A^{(j)} = \I_n + \W \H_j \tra \RBold\inv \H_j\W, \quad j \geq 1,
        \end{equation}
        can then be solved by using an iterative method such as the preconditioned conjugate gradients method (PCG). 

\begin{figure}[htb]
            \centering
            \includegraphics[width=8cm]{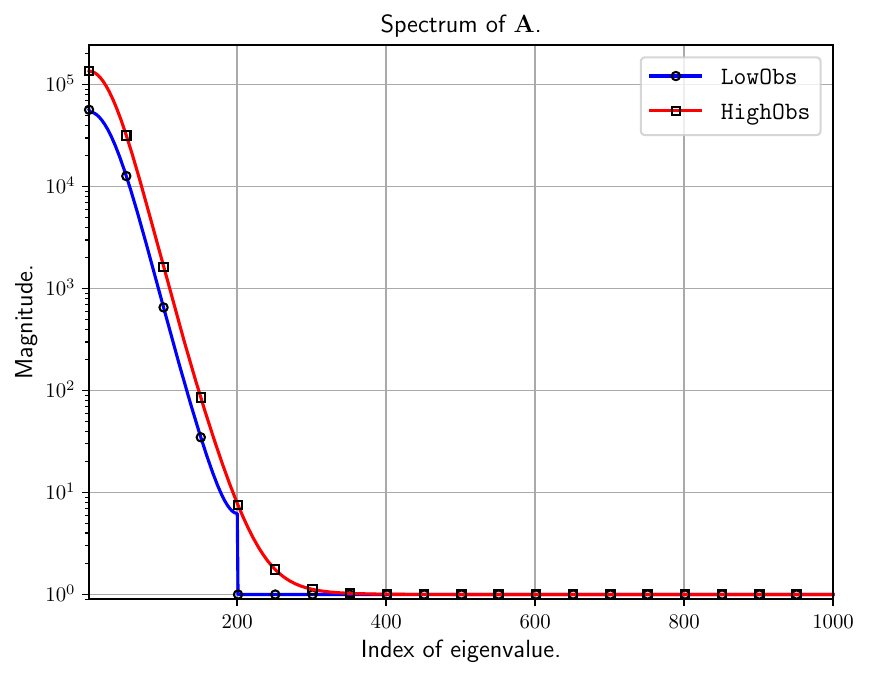}
            \caption{Eigenvalue distribution of $\A$ in the \texttt{LowObs} and \texttt{HighObs} scenarios.}
            \label{fig:spectrum_of_A}
        \end{figure}

        For simplicity, we will use $\A$ referring to $\A^{(1)}$. The error covariance matrix $\B$ is modeled using a finite-difference discretization of the implicit diffusion equation on a regular grid~\cite{GouxGurolEtAl_ConvergenceRDiffusion_2024}, a method widely employed in geophysical applications~\cite{weaver_courtier_2001,WeaverMirouze_Diffusion_2013}. With such a diffusion-modeled covariance $\B$, the matrix vector products with $\W$ can be readily carried out through an operator~\cite{weaver_courtier_2001,GouxGurolEtAl_ConvergenceRDiffusion_2024}. In the diffusion model, we set the the total number of diffusion steps, $M= 6$ and the Daley length-scale $D = 10$, which controls the spatial smoothing properties of the diffusion kernel. $\mathcal{H}$ is chosen as a uniform selection operator ($\mathcal{H}$ is a linear operator),  associated with uniformly distributed observations. The error covariance matrix of observations is set to a diagonal matrix, $\RBold = \sigma_{\RBold} \, \I_m$ with $\sigma_{\RBold} = 10^{-2} > 0$ being the standard deviation. The matrix of interest defined in~\eqref{eq:da_matrix} is a rank $m$ update of the identity. Consequently, the number of observations $m$ has a critical influence on the eigenvalue distribution, i.e. for $n > m$ there are $n-m$ eigenvalues equal to 1. We choose $n= 1000$, and we consider two different scenarios, with a different number of observations: \texttt{LowObs} for $m=200$ ($\mathcal{H} \in \mathbb{R}^{m \times n}$ is a linear selection operator that selects one observation every five grid points), and \texttt{HighObs} for $m=500$ ($\mathcal{H} \in \mathbb{R}^{m \times n}$ is a linear selection operator that selects one observation every two grid points). Figure \ref{fig:spectrum_of_A} shows the spectrum of $\A$ in both cases. Due to its rapidly decaying spectrum, the RSVD is expected to perform well on such problems.

    \subsection{Randomized algorithm to approximate eigenspectrum of an SPD matrix} \label{sec:generalized_SVD_algo}
        In our numerical experiments, we then aim to approximate the eigenspectrum of the symmetric positive definite matrix $\A^{(j)}$. For symmetric matrices, singular and eigenvectors are identical so in this case Algorithm \ref{alg:grsvd}, will actually output an approximate eigenvalue decomposition (EVD). However, the power iteration scheme simplifies for square symmetric matrices, and we rather consider matrices of the form $\ZBold = \A^{q+1} \LBold \G$, with $q \geq 0$. Practically, this reduces to change line 4 in Algorithm \ref{alg:grsvd} to rather compute $\ZBold = \A \QBold$ instead of $\ZBold = \A\tra\A \QBold$. Apart from standard numerical linear algebra operations, the algorithm requires 1 application of $\LBold$ and $q+1$ applications of $\A$ for randomized range-finding (steps 1-5), and one additional application of $\A$ for the low rank approximation computation (steps 6-8). Note that the $m$ dominant eigenvectors of $\A$ are the eigenvectors of the rank $m < n$ matrix $\PhiBold = \W \H_j \tra \RBold\inv \H_j\W$. Consequently, whenever the target rank $k \le m$,
        Algorithm \ref{alg:grsvd} is applied to $\PhiBold$ instead of $\A$.
    
    \subsection{Accuracy of the error bounds} \label{subsec:accuracy_bounds}
        In this section, we evaluate the accuracy of the bounds in Theorem~\ref{th:main} using the VarDA linear system matrix introduced in Section~\ref{sec:3dvar} with \texttt{HighObs} scenario (conclusions in the \texttt{LowObs} scenario are similar). The objective here is twofold: to  numerically illustrate the discussions in Section \ref{sec:application} and to compare the theoretical bounds with the empirical error. For the latest point, the empirical mean of the randomized low-rank approximation error is computed out of 20 independent runs of Algorithm~\ref{alg:grsvd}.
        
        For comparison, we plot the following normalized quantity
        \begin{equation*}
            \dfrac{\norm{[\I_m - \Proj(\ZBold)] \A}_F }{ \norm{\overline{\SigmaBold}_k}_F} - 1,
        \end{equation*}
        which approaches zero as the low-rank approximation error approaches the optimal value $\norm{\overline{\SigmaBold}_k}_F$.
        
\begin{figure}[ht!]
            \centering
            \begin{subfigure}{0.495\textwidth}
                \centering
                \includegraphics[width=7.15cm]{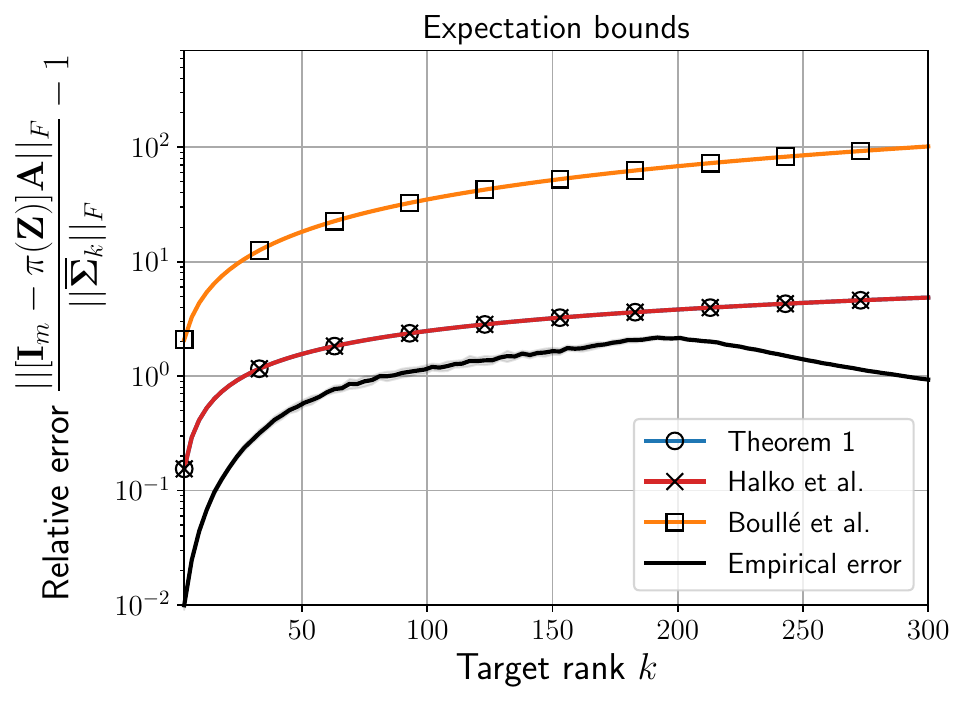}
                \caption{Expectation bounds.}
                \label{fig:rsvd_vs_k:expect}
            \end{subfigure}
            \begin{subfigure}{0.495\textwidth}
                \centering
                \includegraphics[width=7.15cm]{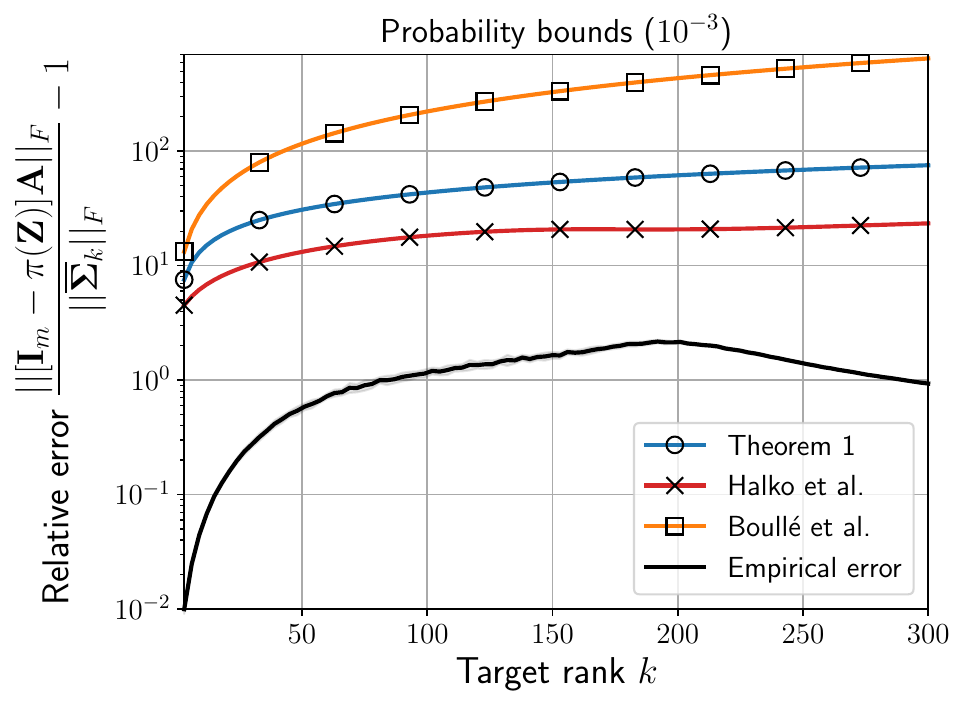}
                \caption{Probability bounds.}
                \label{fig:rsvd_vs_k:proba}
            \end{subfigure}
            \caption{Bounds for the low-rank approximation error versus the target rank $k$ with $p=10$.}
            \label{fig:rsvd_vs_k}
        \end{figure}

\begin{figure}[ht!]
            \centering
            \begin{subfigure}{0.495\textwidth}
                \centering
                \includegraphics[width=7.15cm]{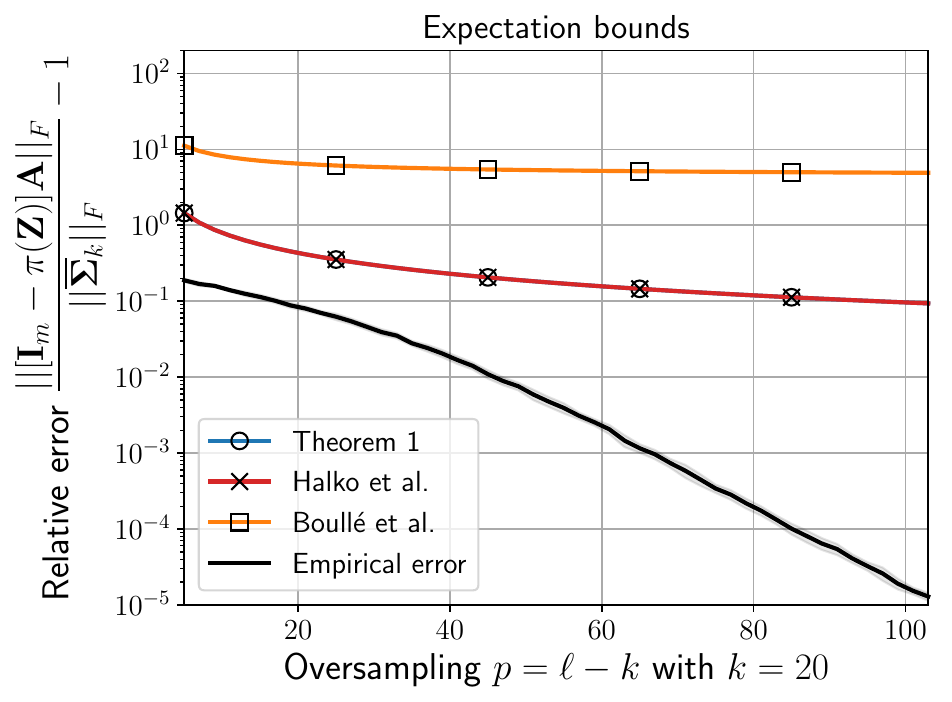}
                \caption{Expectation bounds.}
                \label{fig:rsvd_vs_p:expect}
            \end{subfigure}
            \begin{subfigure}{0.495\textwidth}
                \centering
                \includegraphics[width=7.15cm]{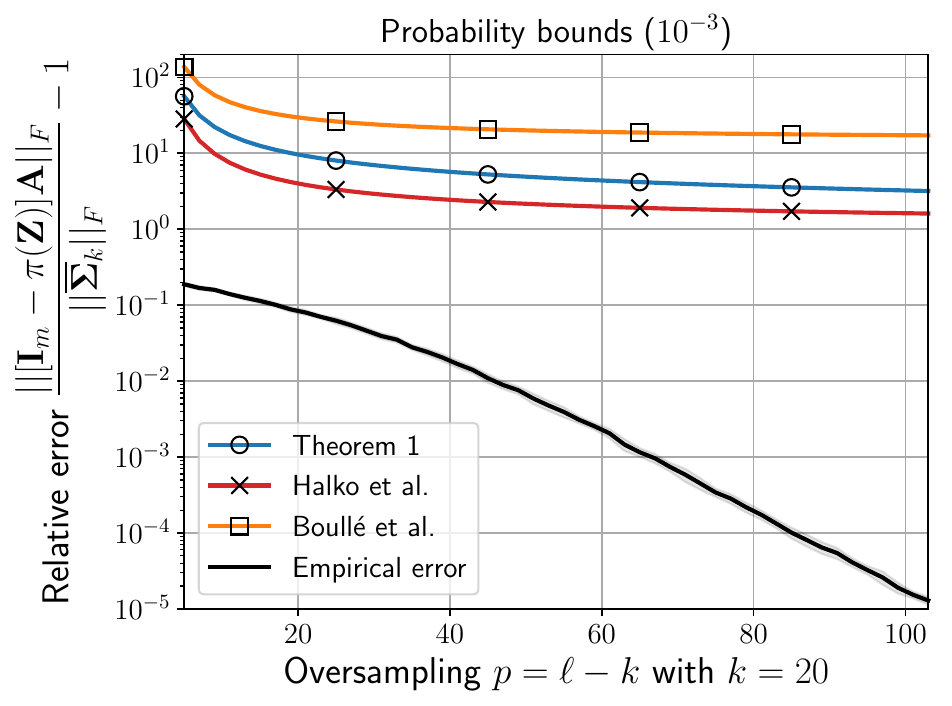}
                \caption{Probability bounds.}
                \label{fig:rsvd_vs_p:proba}
            \end{subfigure}
            \caption{Bounds for the low-rank approximation error versus the oversampling parameter $p = \ell - k$ with $k=20$.}
            \label{fig:rsvd_vs_p}
        \end{figure}

        Let us first analyze the bounds accuracy in the RSVD case, that is $\LBold = \I_n$ and $q=0$. Figures \ref{fig:rsvd_vs_k} and \ref{fig:rsvd_vs_p} show the comparison of the bounds in expectation and in probability, respectively versus the target rank $k$ and the oversampling $\ell - k$. As pointed out in Section \ref{sec:application}, expectation bounds of Theorem \ref{th:main} and~\cite{HalkoMartinssonEtAl_FindingStructureRandomness_2011} are identical while our probability bounds is slightly worse. The accuracy loss from the generalization in~\cite{BoulleTownsend_GeneralizationRandomizedSingular_2022} is also clearly visible. Overall, given the small empirical standard deviation observed (gray halo around the empirical error), all the probability bounds tend to severely overestimate the actual error. Another important comment is that the behavior of the bounds with respect to the oversampling does not match the one of the actual error. Indeed, the bounds predict a decrease in $1/\sqrt{\ell-k}$ while we see in Figure \ref{fig:rsvd_vs_p} that the actual decrease seems exponential (linear in semi-log scale).
    
\begin{figure}[ht!]
            \centering
            \begin{subfigure}{0.495\textwidth}
                \centering
                \includegraphics[width=7.15cm]{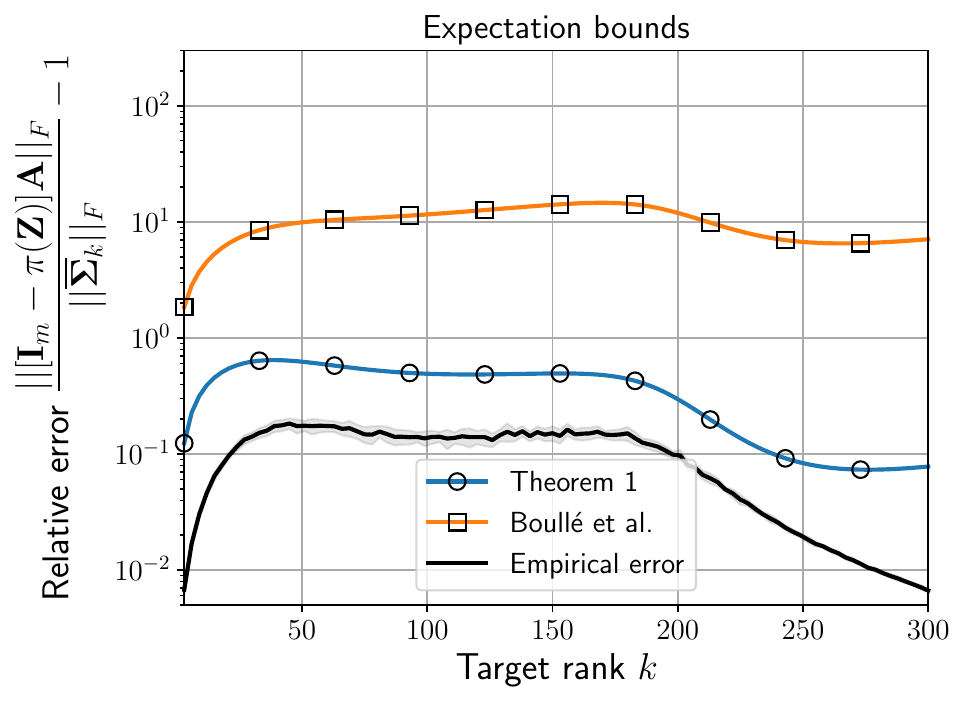}
                \caption{Expectation bounds.}
                \label{fig:grsvd_vs_k:expect}
            \end{subfigure}
            \begin{subfigure}{0.495\textwidth}
                \centering
                \includegraphics[width=7.15cm]{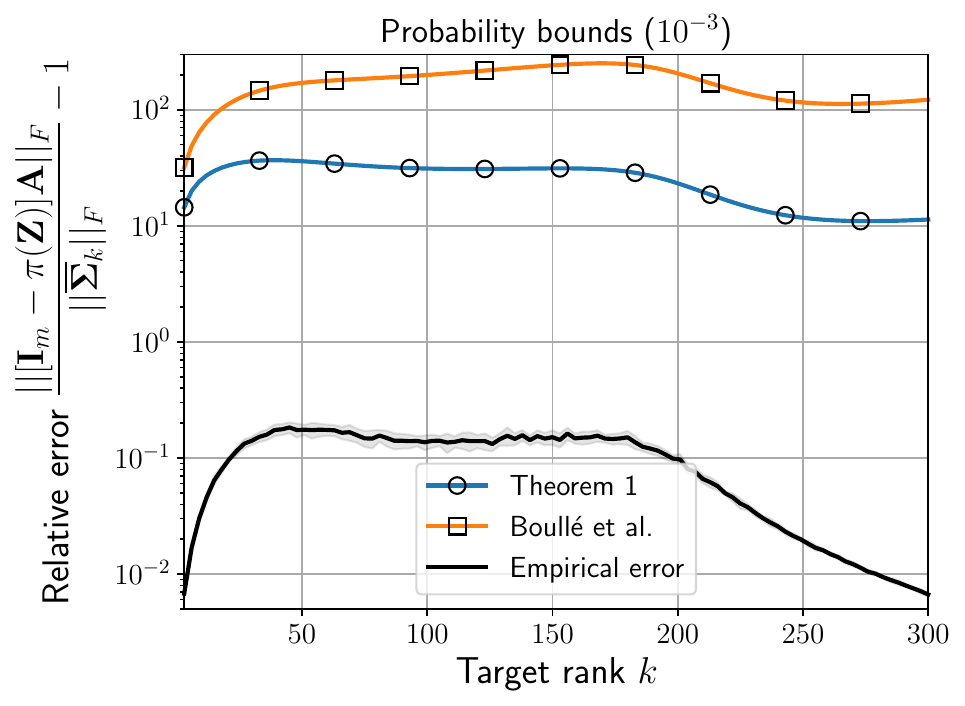}
                \caption{Probability bounds.}
                \label{fig:grsvd_vs_k:proba}
            \end{subfigure}
            \caption{Bounds for the low-rank approximation error versus the target rank $k$ with $p=10$ and $\LBold = \W$.}
            \label{fig:grsvd_vs_k}
        \end{figure}

\begin{figure}[ht!]
            \centering
            \begin{subfigure}{0.495\textwidth}
                \centering
                \includegraphics[width=7.15cm]{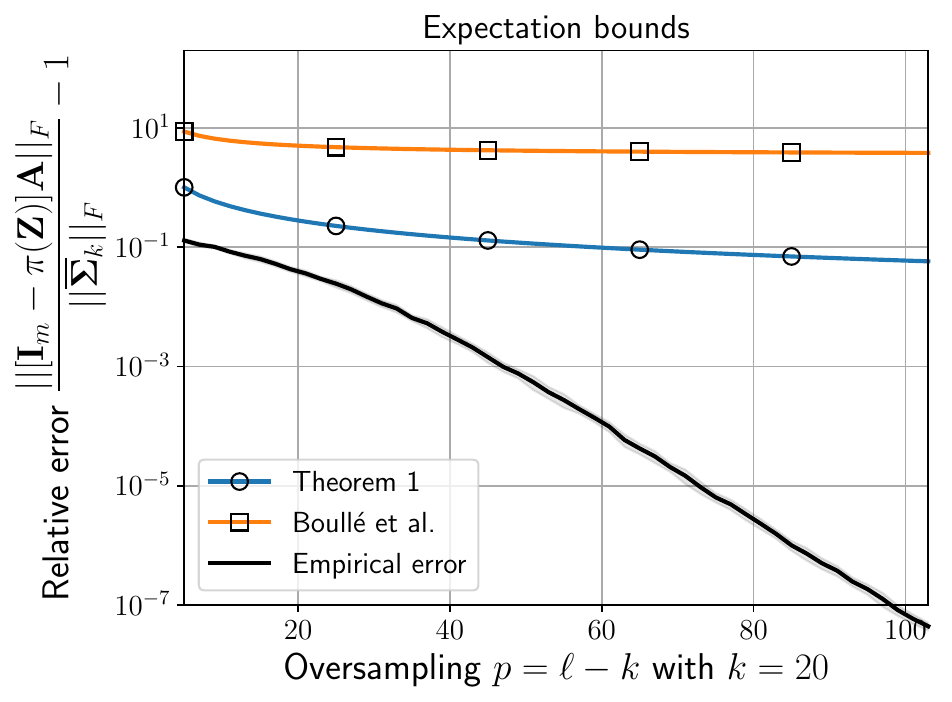}
                \caption{Expectation bounds.}
                \label{fig:grsvd_vs_p:expect}
            \end{subfigure}
            \begin{subfigure}{0.495\textwidth}
                \centering
                \includegraphics[width=7.15cm]{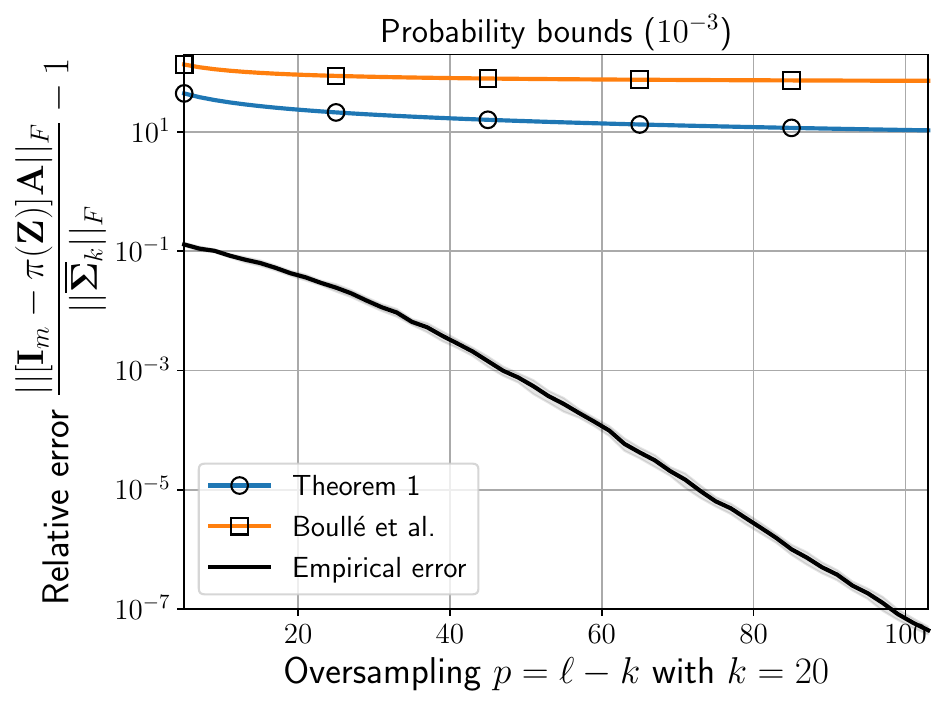}
                \caption{Probability bounds.}
                \label{fig:grsvd_vs_p:proba}
            \end{subfigure}
            \caption{Bounds for the low-rank approximation error versus the oversampling parameter $p = \ell - k$ with $k=20$ and $\LBold = \W$.}
            \label{fig:grsvd_vs_p}
        \end{figure}

        Secondly, we propose a comparison in the generalized SVD setting, where $\LBold = \W$. The relevance of this particular choice for $\LBold$ is discussed in the subsequent section. In this case, the bounds from~\cite{HalkoMartinssonEtAl_FindingStructureRandomness_2011} are not longer usable, and we provide a comparison solely with the bounds in~\cite{BoulleTownsend_GeneralizationRandomizedSingular_2022}. In Figure \ref{fig:grsvd_vs_k}, we see that the expectation bound captures well the trend of the empirical error. By contrast, and similarly to the RSVD case above, Figure \ref{fig:grsvd_vs_p} illustrates that the bounds still fail to capture the behavior with respect to the oversampling. The gap between our bound and the ones from Boullé et al. increases with $k$ as highlighted in Section \ref{sec:application:grsvd}. Again, the probability bounds predict low rank approximation errors that are significantly overestimated.
        
    \subsection{Different choices of the matrix \textbf{L}} \label{subsec:choices_for_L}
    
        In this section, we investigate different choices of $\LBold$. Specifically, we aim to answer the following questions: 
        \begin{itemize}
            \item Performing one exact subspace iteration can be done by considering either $q = 1$ and $\mathbf{L} = \mathbf{I}_n$ or equivalently $q = 0$ and $\mathbf{L} = \PhiBold$. Since applying $\PhiBold$ is costly, the question is: can we choose a computationally efficient choice of $\mathbf{L}$ to reduce the computational cost of Algorithm~\ref{alg:grsvd} while maintaining similar accuracy? This is the object of Section \ref{subsec:approximate_subspace_iteration}.
            \item Can we achieve better accuracy with a different choice of $\mathbf{L}$ compared to the standard setting ($q = 0$ and $\mathbf{L} = \mathbf{I}_n$) without significantly increasing the computational cost? This is the object of Sections \ref{subsec:available_info} and \ref{subsec:combination}.
        \end{itemize}
    
        \subsubsection{Using   \textbf{W} in the subspace iteration:} \label{subsec:approximate_subspace_iteration}
        
            The different settings of Algorithm \ref{alg:grsvd} we consider are gathered in Table \ref{tab:algo_settings}. The first setting is the standard RSVD, which serves as a reference. The second and third settings illustrate the approximate subspace iterations where $\PhiBold$ is replaced by $\W$. Note that the range of $\PhiBold$ is in the range of $\W$. Consequently, $\W$ satisfies the properties highlighted for $\LBold$ to result in efficient low rank approximation algorithms. Finally, the last setting is the exact subspace iteration with two steps, which also serves as a reference.
    
            \begin{table}[ht!]
                \centering
                \begin{tabular}{c||c|c|c|c}
                    Matrix $\LBold$ & $\I_n$ & $\W$ & $\W^2 = \B$ & $\I_n$ \\ \hline
                    Power iteration steps $q$ & 0 & 0 & 0 & 1
                \end{tabular}
                \caption{Settings of Algorithm \ref{alg:grsvd} for different approaches.}
                \label{tab:algo_settings}
            \end{table}
    
\begin{figure}[ht!]
                \centering
                \begin{subfigure}{0.495\textwidth}
                    \centering
                    \includegraphics[width=7.15cm]{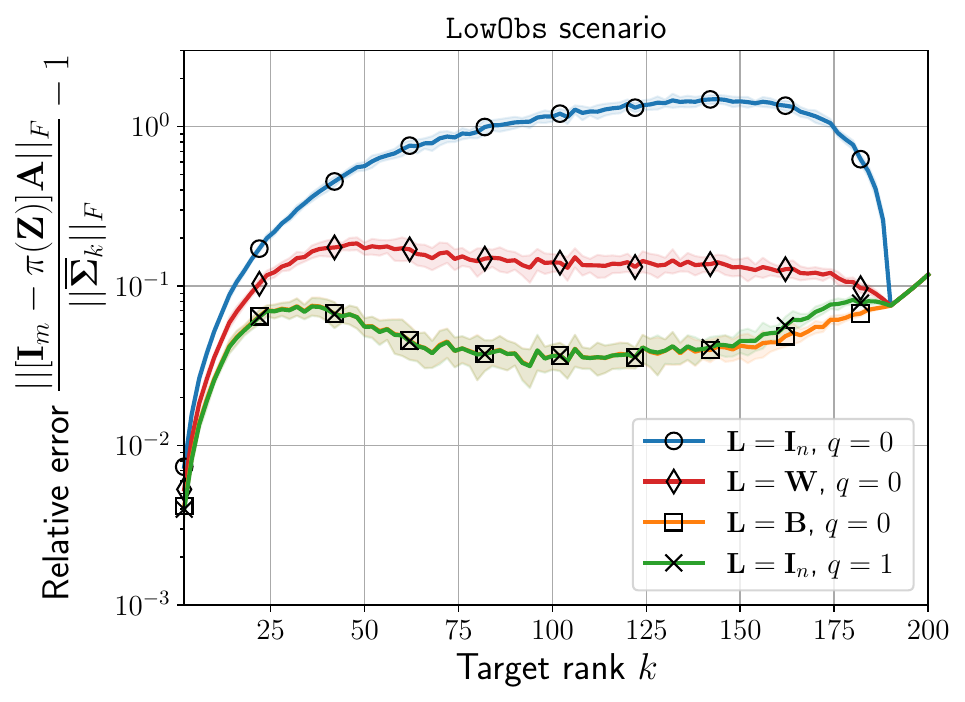}
                    \caption{\texttt{LowObs}: $m=200$.}
                    \label{fig:approx_subspace_iteration_vs_k:low}
                \end{subfigure}
                \begin{subfigure}{0.495\textwidth}
                    \centering
                    \includegraphics[width=7.15cm]{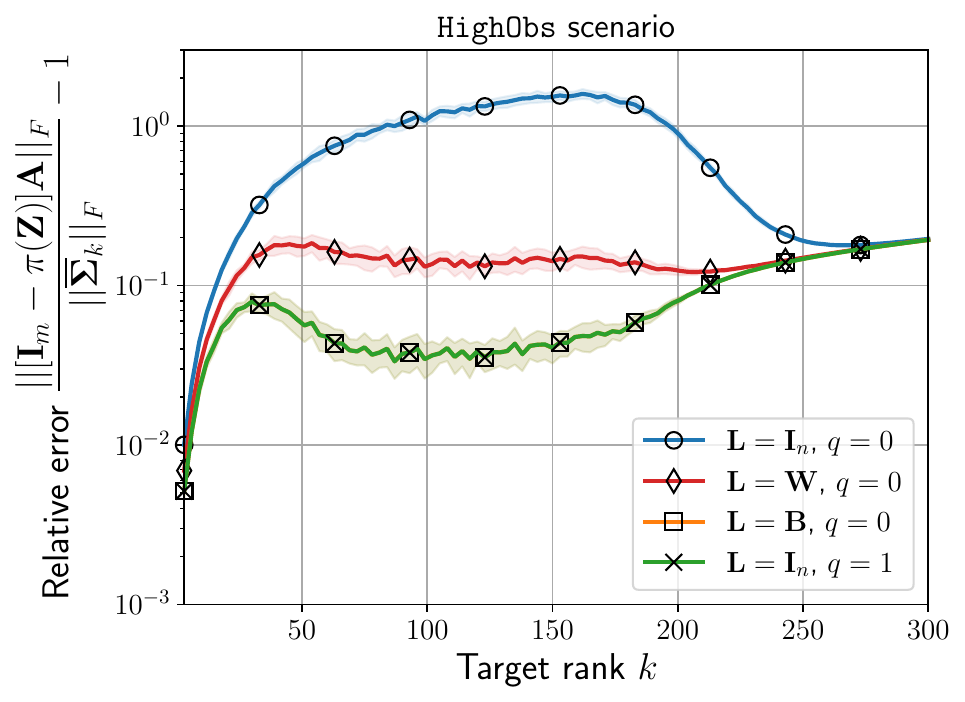}
                    \caption{\texttt{HighObs}: $m=500$.}
                    \label{fig:approx_subspace_iteration_vs_k:high}
                \end{subfigure}            
                \caption{Empirical low-rank approximation error versus the target rank with $p=10$ and different Algorithm \ref{alg:grsvd} settings.}
                \label{fig:approx_subspace_iteration_vs_k}
            \end{figure}

\begin{figure}[ht!]
                \centering
                \begin{subfigure}{0.495\textwidth}
                    \centering
                    \includegraphics[width=7.15cm]{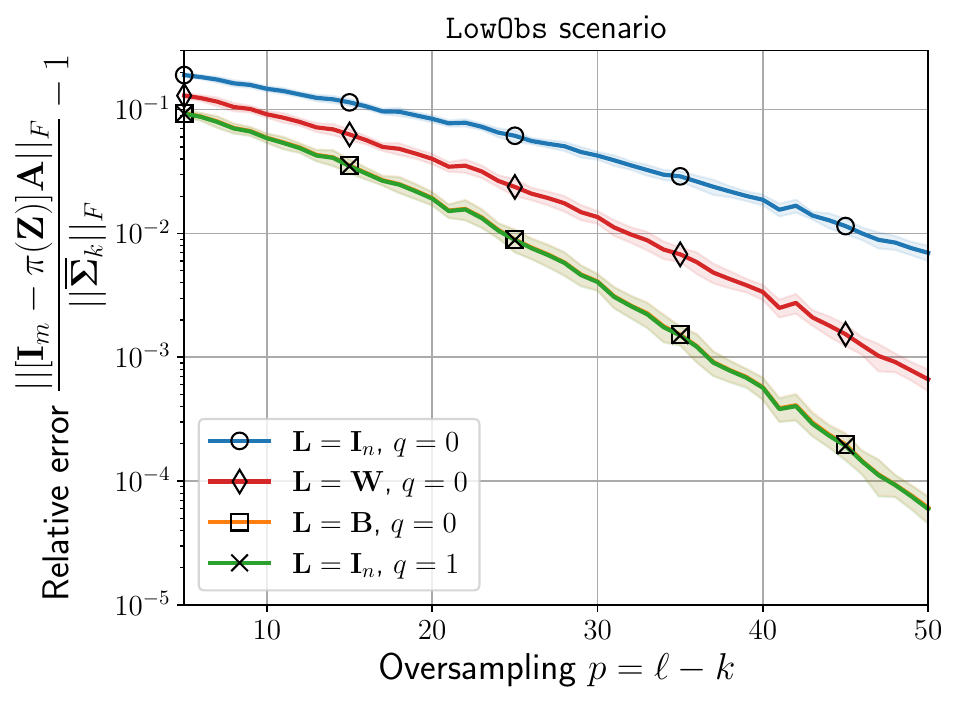}
                    \caption{\texttt{LowObs}: $m=200$.}
                    \label{fig:approx_subspace_iteration_vs_p:low}
                \end{subfigure}
                \begin{subfigure}{0.495\textwidth}
                    \centering
                    \includegraphics[width=7.15cm]{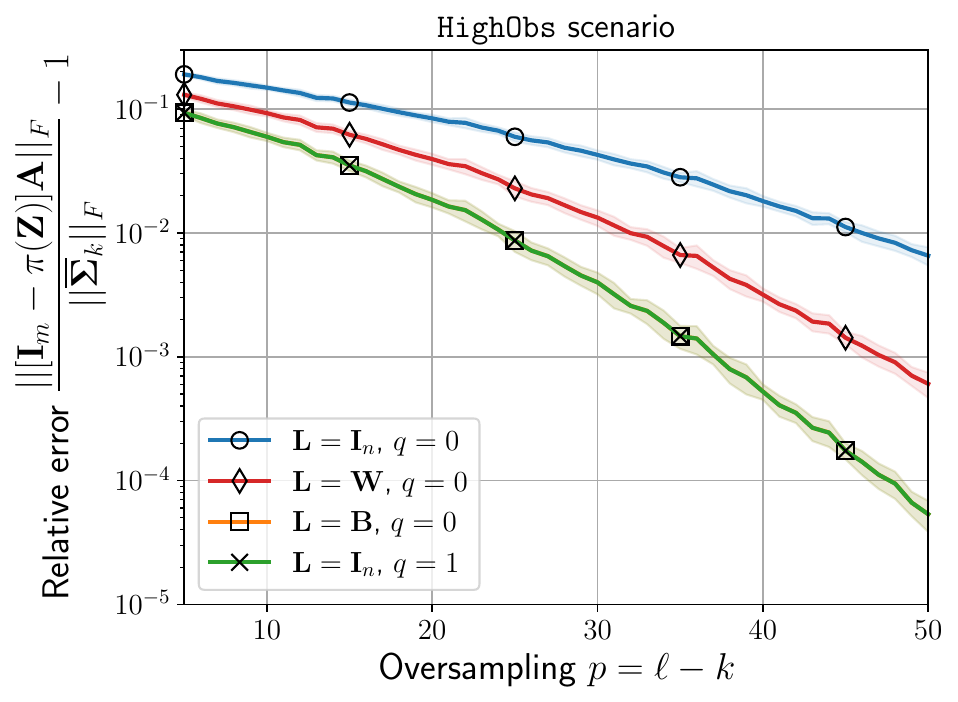}
                    \caption{\texttt{HighObs}: $m=500$.}
                    \label{fig:approx_subspace_iteration_vs_p:high}
                \end{subfigure}            
                \caption{Empirical low-rank approximation error versus the oversampling with $k=20$ and different Algorithm \ref{alg:grsvd} settings.}
                \label{fig:approx_subspace_iteration_vs_p}
            \end{figure}

            The resulting low rank approximation errors in both the \texttt{LowObs} and \texttt{HighObs} scenarios are shown in Figures \ref{fig:approx_subspace_iteration_vs_k} and \ref{fig:approx_subspace_iteration_vs_p}, versus the target rank $k$ and the oversampling $\ell - k$ respectively. First, it is clear that using powers of $\W$ for $\LBold$ improves the resulting low rank approximation. Furthermore, in both scenarios, $\LBold = \W^2 = \B$ actually yields similar error as using $\A$, while being computationally cheaper. This means that for this particular problem, the dominant eigenmodes of $\A$ are mostly carried by $\B$. Figure \ref{fig:approx_subspace_iteration_vs_k} shows that when the target rank $k$ approaches the number of observations $m$, the different settings tend to provide similar errors. This means that up to a certain point, the additional information brought by $\LBold$ is no longer relevant. This is consistent with the spectrum of $\A$ in Figure \ref{fig:spectrum_of_A}: when $k$ gets close to $m$, the eigenmodes that we try to approximate are associated to eigenvalues close to 1, and are thus hard to isolate from the remaining $n-m$ eigenmodes associated to 1. This fact also explains why the low rank approximation error increases with $k$. In Figure \ref{fig:approx_subspace_iteration_vs_p}, we observe that the oversampling allows to significantly improve the performance. The hierarchy between the settings is similar, but the slope becomes sharper with the powers of $\W$. Consequently, the oversampling is more beneficial, relatively, when additional information is brought via $\LBold$. Finally, we also would like to point out that the empirical standard deviation observed is small, implying that the results of Algorithm \ref{alg:grsvd} are fairly robust.
        
        \subsubsection{Using available approximations:} \label{subsec:available_info}
    
           Let us assume that there exist already available approximations $\widehat{\V}_k \in \Rnk$ and $\widehat{\SigmaBold}_k \in \Rkk$ of $\V_k$ and $\SigmaBold_k$ respectively, where $\widehat{\SigmaBold}_k = \diag(\hat{\sigma}_1, \dots, \hat{\sigma}_k)$ and $\hat{\sigma}_1 \geq \dots \geq\hat{\sigma}_k > 0$. This follows the discussion in Section \ref{sec:application:available_approx}, where we considered, for $\beta > 0$, matrices of the form:
            \begin{equation*}
                \LBold_\beta =
                \widehat{\V}_k \widehat{\SigmaBold}_k \widehat{\V}_k\tra +
                \beta \hat{\sigma}_k (\I_n - \widehat{\V}_k \widehat{\V}_k\tra)
                = \widehat{\V}_k \left(\widehat{\SigmaBold}_k - \beta \hat{\sigma}_k \I_k\right) \widehat{\V}_k\tra + \beta \hat{\sigma}_k \I_n.
            \end{equation*}
            The last form is more efficient computationally. Indeed the matrix-vector products with $\LBold_\beta$ costs approximately $(4k+1)n$ operations and could be performed efficiently even for very large scale problems. In DA, these approximations could come from previous systems solutions in the sequence of linear systems \eqref{eq:da_linear_system}. For instance, if system $j$ is solved using a Krylov method, then approximate eigenpairs can be computed from the Krylov subspace~\cite{Saad_IterativeMethodsLinear_2019, FishNoceTremWrig09}. In operational settings, such information is readily available at negligible cost and is, in fact, commonly used for preconditioning~\cite{FishNoceTremWrig09}. In the following experiments, we obtain $\widehat{\V}_k$ and $\widehat{\SigmaBold}_k$ by using the RSVD applied to the SPD matrix $\PhiBold$, with an oversampling parameter set to $5$. 
   
            Let us begin by studying the following settings:
            \begin{itemize}
                \item $\LBold_0$: as discussed in Section \ref{sec:application:available_approx}, this yields a deterministic algorithm where the standard power iteration is performed with initial guess $\widehat{\V}_k$. This choice will serve as a reference, and our objective is to study whether maintaining randomization (i.e. $\beta > 0$) can yield better performance.
                \item $\LBold_{0.01}$: here we allow an exploration of the orthogonal of $\range(\widehat{\V}_k)$ while being relatively confident in the quality of $\widehat{\V}_k$ as an approximation of $\V_k$.
                \item $\LBold_{1}$: this choice ponders negatively the confidence in the quality of $\widehat{\V}_k$, and more weight is given to the orthogonal of $\range(\widehat{\V}_k)$.
            \end{itemize}
    
\begin{figure}[ht!]
                \centering
                \begin{subfigure}{0.495\textwidth}
                    \centering
                    \includegraphics[width=7.15cm]{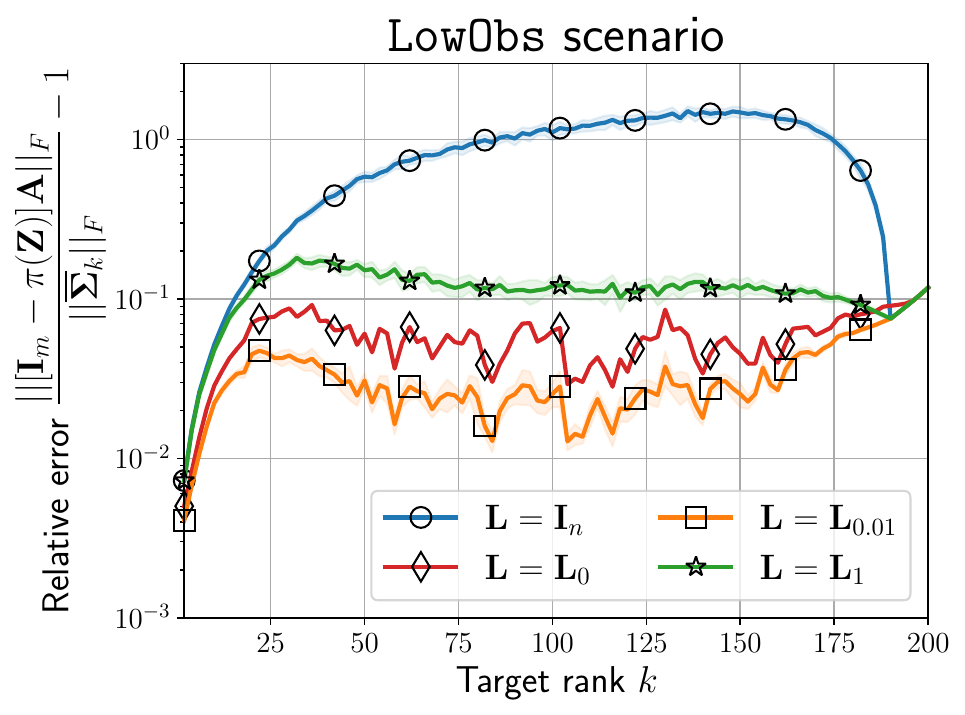}
                    \caption{\texttt{LowObs}: $m=200$.}
                    \label{fig:available_approx_vs_k:low}
                \end{subfigure}
                \begin{subfigure}{0.495\textwidth}
                    \centering
                    \includegraphics[width=7.15cm]{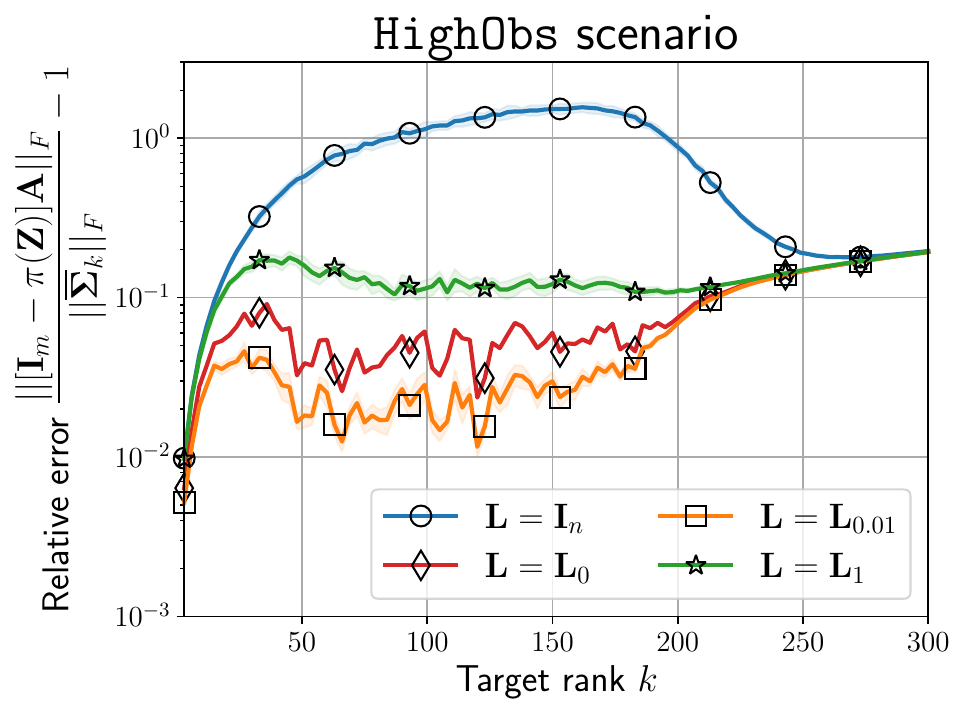}
                    \caption{\texttt{HighObs}: $m=500$.}
                    \label{fig:available_approx_vs_k:high}
                \end{subfigure}            
                \caption{Empirical low-rank approximation error versus the target rank with $p=10$ and different Algorithm \ref{alg:grsvd} settings.}
                \label{fig:available_approx_vs_k}
            \end{figure}

\begin{figure}[ht!]
                \centering
                \begin{subfigure}{0.495\textwidth}
                    \centering
                    \includegraphics[width=7.15cm]{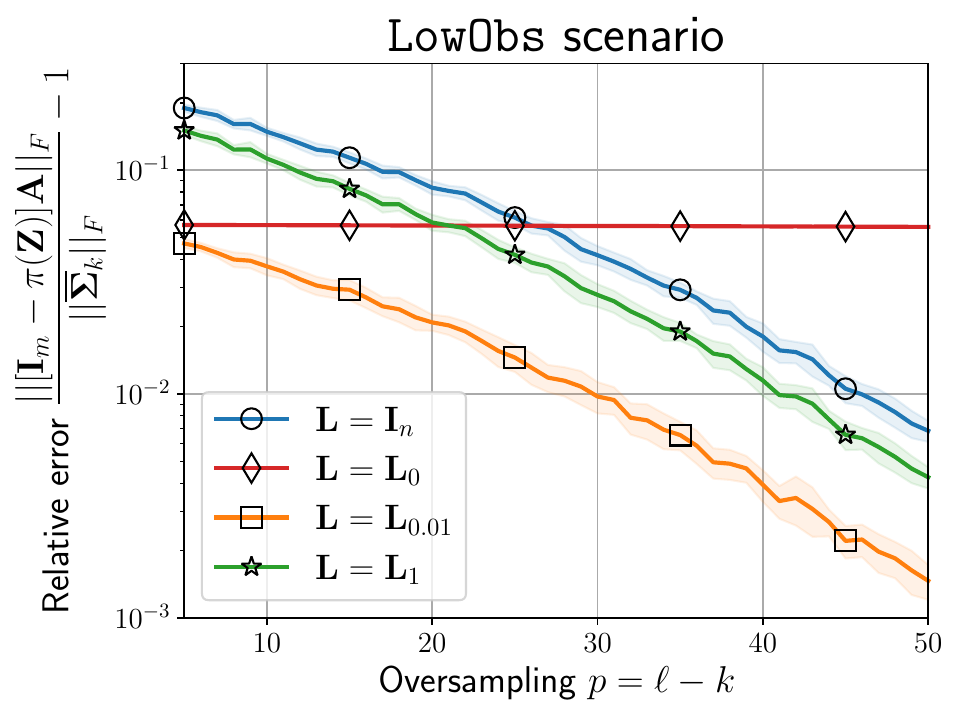}
                    \caption{\texttt{LowObs}: $m=200$.}
                    \label{fig:available_approx_vs_p:low}
                \end{subfigure}
                \begin{subfigure}{0.495\textwidth}
                    \centering
                    \includegraphics[width=7.15cm]{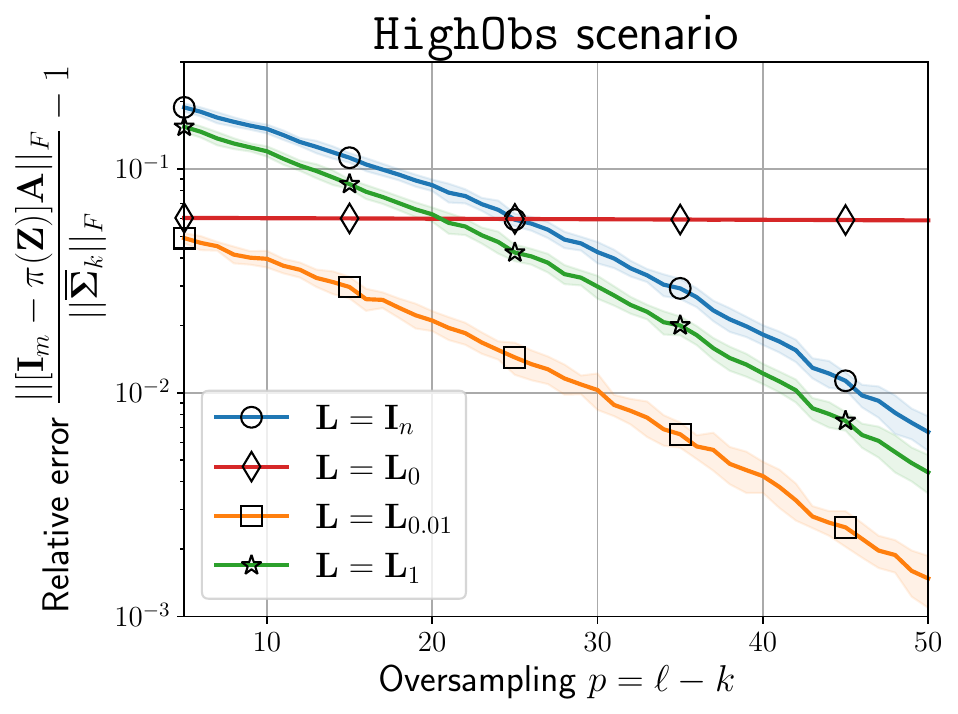}
                    \caption{\texttt{HighObs}: $m=500$.}
                    \label{fig:available_approx_vs_p:high}
                \end{subfigure}            
                \caption{Empirical low-rank approximation error versus the oversampling with $k=20$ and different Algorithm \ref{alg:grsvd} settings.}
                \label{fig:available_approx_vs_p}
            \end{figure}

            The results are shown in Figures \ref{fig:available_approx_vs_k} and \ref{fig:available_approx_vs_p}, versus the target rank $k$ and the oversampling $\ell - k$ respectively. We observe that using $\LBold_{0.01}$ yields better results than the reference $\LBold_0$. The fact that $\LBold_{1}$ performs worse shows indirectly that the approximations used to construct $\LBold_{\beta}$ are accurate, and that smaller values of $\beta$ must be preferred. Also, all the configurations yield better results than $\LBold = \I_n$, showing that using $\LBold_{\beta}$ brings relevant additional information. In Figure \ref{fig:available_approx_vs_p}, $\LBold_0$ yields constant results, i.e. $\range(\LBold_0 \G)$ is constant for any value of $\ell$ as expected. By contrast with the observation in Figure \ref{fig:approx_subspace_iteration_vs_p}, the decrease in low rank approximation error with the oversampling has the same slope regardless of the choice for $\beta > 0$. The gain with increasing values of the oversampling of $\LBold_{0.01}$ compared to $\LBold_0$ is of great interest. Indeed, in practice the number of eigenpairs available is not user defined, so the use of $\LBold_\beta$ with $\beta > 0$ combined with an oversampling may significantly improve the result compared to simply using $\widehat{\V}_k$.
    
        \subsubsection{Combining  approximate information and \textbf{W}:} \label{subsec:combination}
           As a last experiment,  we propose to combine the use of $\W$ and $\LBold_\beta$. In fact, we have seen that each approach allows improving the performance while involving an extra numerical cost which is affordable. We thus consider the following settings:
            \begin{itemize}
                \item $\LBold := \W$: here as a reference for the use of $\W$ alone, similar to Section \ref{subsec:approximate_subspace_iteration}.
                \item $\LBold := \LBold_{0.01}$: also used as a reference to have the error with $\LBold_\beta$ alone, similar to Section \ref{subsec:available_info}.  We have selected $\beta = 0.01$ as it yielded the better performance.
                \item $\LBold := \W \LBold_{0.01}$: the case of interest where the information from both $\W$ and $\LBold_\beta$ is combined.
            \end{itemize}
    
\begin{figure}[ht!]
                \centering
                \begin{subfigure}{0.495\textwidth}
                    \centering
                    \includegraphics[width=7.15cm]{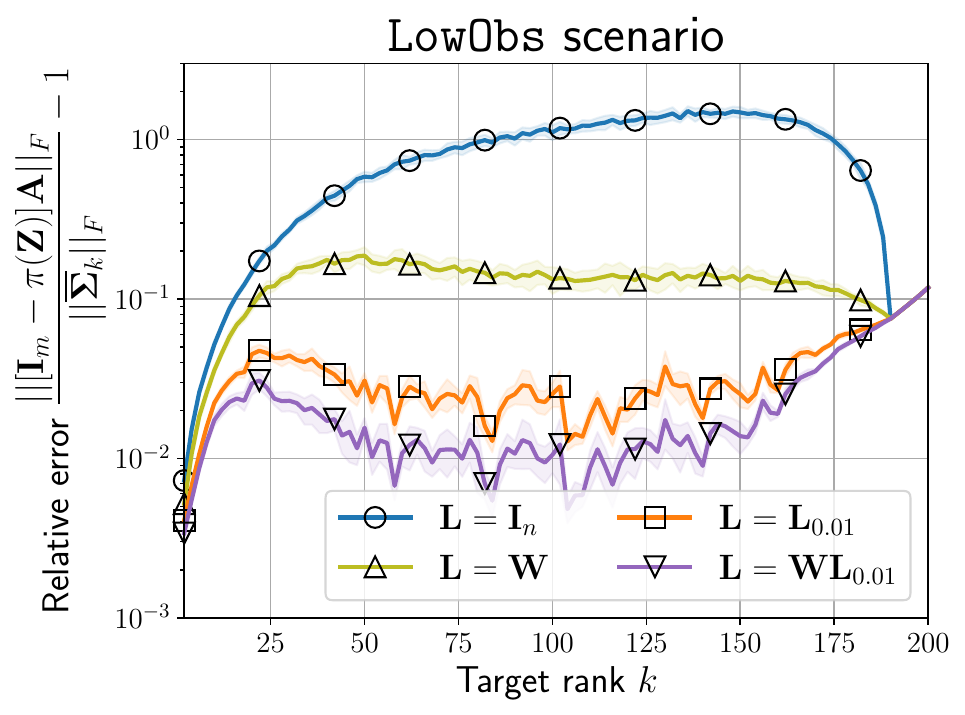}
                    \caption{\texttt{LowObs}: $m=200$.}
                    \label{fig:b_available_approx_vs_k:low}
                \end{subfigure}
                \begin{subfigure}{0.495\textwidth}
                    \centering
                    \includegraphics[width=7.15cm]{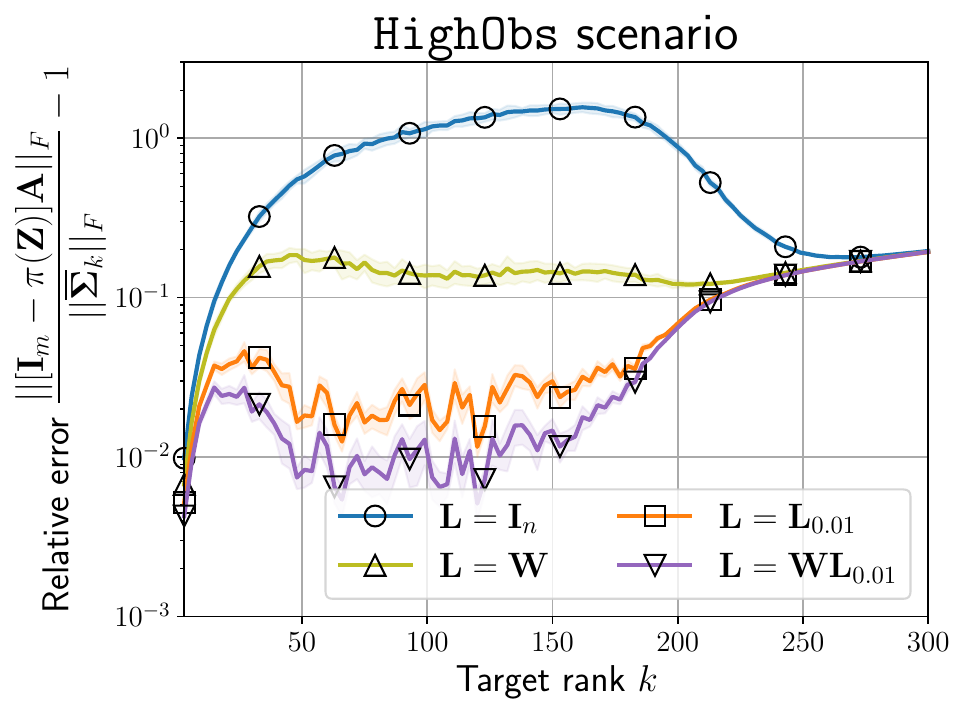}
                    \caption{\texttt{HighObs}: $m=500$.}
                    \label{fig:b_available_approx_vs_k:high}
                \end{subfigure}            
                \caption{Empirical low-rank approximation error versus the target rank with $p=10$ and different Algorithm \ref{alg:grsvd} settings.}
                \label{fig:b_available_approx_vs_k}
            \end{figure}

\begin{figure}[ht!]
                \centering
                \begin{subfigure}{0.495\textwidth}
                    \centering
                    \includegraphics[width=7.15cm]{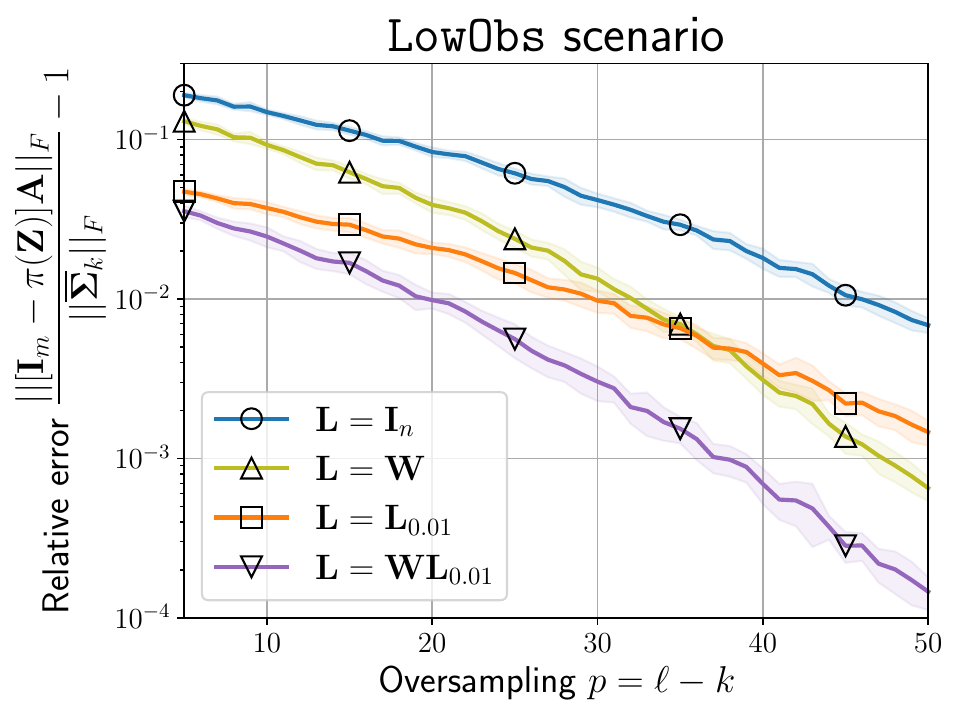}
                    \caption{\texttt{LowObs}: $m=200$.}
                    \label{fig:b_available_approx_vs_p:low}
                \end{subfigure}
                \begin{subfigure}{0.495\textwidth}
                    \centering
                    \includegraphics[width=7.15cm]{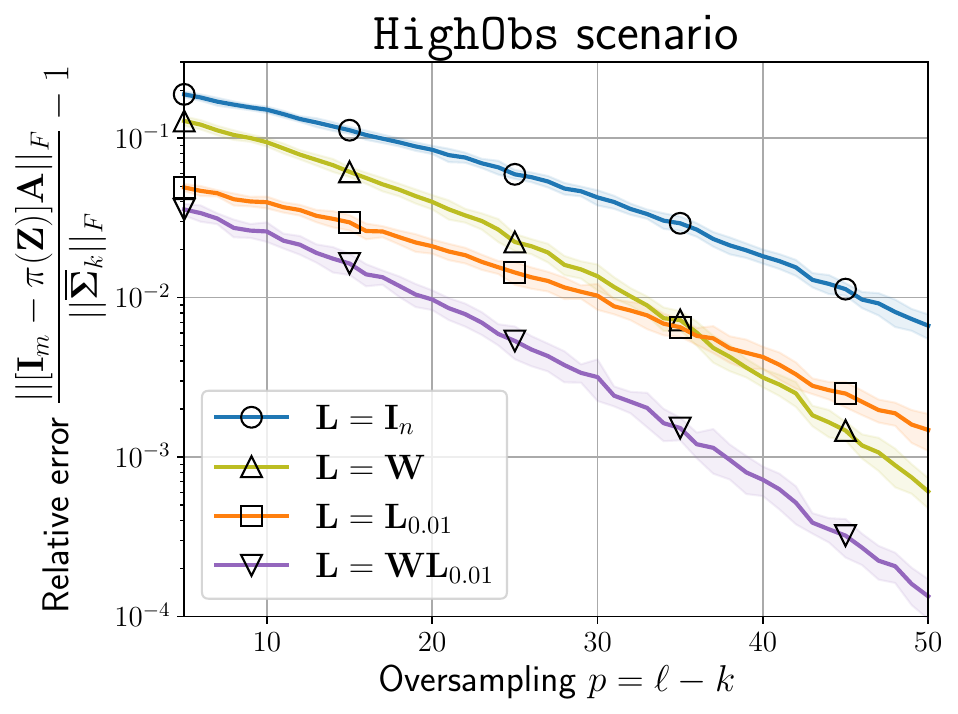}
                    \caption{\texttt{HighObs}: $m=500$.}
                    \label{fig:b_available_approx_vs_p:high}
                \end{subfigure}            
                \caption{Empirical low-rank approximation error versus the oversampling with $k=20$ and different Algorithm \ref{alg:grsvd} settings.}
                \label{fig:b_available_approx_vs_p}
            \end{figure}

            The results are shown in Figures \ref{fig:b_available_approx_vs_k} and \ref{fig:b_available_approx_vs_p}. The fact that the combination of both $\W$ and $\LBold_\beta$ yields better results than each of them separately illustrates that the two approaches bring different eigen information on $\A$. Again, in operational contexts, applying $\W$ or $\LBold_\beta$ is computationally less expensive than applying $\A$ and is thus numerically very appealing. We see in Figure \ref{fig:b_available_approx_vs_p} the difference of slope observed above. The use of $\W$ results in a sharper slope, which increases the gain in low rank approximation error obtained with oversampling.
    
    As a conclusion, the usage of non-standard choices of $\LBold$ with $q=0$ has a strong practical potential. In VarDA, these choices can rely on either the error covariance matrix $\B$, or approximate eigen information computed from the solution of previous systems. In any case, such construction requires negligible extra computational cost but can lead to a great increase in the approximate eigen information obtained from Algorithm \ref{alg:grsvd}, either used separately as in Sections \ref{subsec:approximate_subspace_iteration} and \ref{subsec:available_info} or combined in Section \ref{subsec:combination}. We emphasize that the numerical experiments presented in this section are problem-dependent (including the choice of $\LBold$); however, they clearly demonstrate the potential of using non-standard Gaussian distributions within randomized algorithms and provide insights on the choice of $\LBold$ for other applications.
    
\section{Conclusions} \label{Sec:conclusions}
    We have proposed a general theoretical framework for the stochastic analysis of the low-rank approximation error where we have relaxed the assumptions on the covariance matrix. The proposed bounds in expectation and in probability come with a clear interpretations allowing to identify the key elements resulting in efficient randomized low rank approximation error. More importantly, the provided generalization allows to recover similar bounds from the literature, thus providing a versatile tool for future randomized algorithms analysis. In particular, we recover the expectation bound for the RSVD from \cite{HalkoMartinssonEtAl_FindingStructureRandomness_2011}, and actually get less pessimistic bounds in the case of the randomized power iteration \cite{Gu_SubspaceIterationRandomization_2015} and generalized RSVD \cite{BoulleTownsend_GeneralizationRandomizedSingular_2022}.
    Finally, we have illustrated our result on a data assimilation problem. First, we have studied the accuracy of the bounds. Then, we have proposed numerical experiments in two different situations, using either an approximation of $\A$ or an available approximation of $\V_k$. In both cases, we have identified covariance matrices that enable us to improve the overall performance. This highlights that using the structure of the problem to design the covariance matrix does improve the performance of randomized low-rank approximation methods. It is clear that a similar generalized error analysis can be done for the 2-norm, but also for other forms of low rank approximations such as the Nyström method for symmetric positive semi-definite matrices. 
    
    Future work will focus on incorporating these algorithms into a preconditioner for the operational variational data assimilation framework within the Object-Oriented Prediction System (OOPS) used by European Centre for Medium Weather Forecasts and M\'et\'eo-France.

\section*{Acknowledgements}

The authors would like to thank the two anonymous reviewers for their constructive comments and suggestions, which have helped improve the clarity and quality of this manuscript. 
\appendix

\section{Proof of Theorem~\ref{th:main:general}} \label{app:proofs} 
    
    Let us begin with Lemma \ref{lem:1} which is a generalization of the analysis in terms of principal angles proposed in~\cite{Saibaba_RandomizedSubspaceIteration_2019}. For readability, let us recall the partitioning of the SVD of $\A$ presented in \eqref{eq:partitioning_svd}, which will be intensively used in the proofs. For a given integer $k \leq \min\set{m, n}$, we consider the following ,
    \begin{equation*}
        \A =
        \begin{bmatrix}
            \U_k & \overline{\U}_k
        \end{bmatrix}
        \begin{bmatrix}
            \SigmaBold_k & \\
            & \overline{\SigmaBold}_k
        \end{bmatrix}
        \begin{bmatrix}
            \V_k & \overline{\V}_k
        \end{bmatrix}\tra,
    \end{equation*}
    
    \begin{lemma}  \label{lem:1}
        Let $\A \in \Rmn$ be an arbitrary matrix and $\ZBold \in \Rml$ a full column rank matrix such that $\ell \leq \rank(\A)$. Let us denote $\ZBold_k = \U_k\tra \ZBold \in \R^{k \times \ell}$ and $\overline{\ZBold}_k = \overline{\U}_k\tra \ZBold \in \R^{(m-k) \times \ell}$.

        \vspace{0.25cm}

        \noindent \noindent For any integer $k \leq \ell$, if $\rank(\ZBold_k) = k$, then one has
        \begin{equation} \label{eq:sine_bound_yd}
            \U_k\tra [\I_n - \Proj(\ZBold)] \U_k ~ \lleq ~ \T_k\tra \T_k,
        \end{equation}
        where $\T_k = \overline{\ZBold}_k \ZBold_k\pinv \in \R^{(m-k) \times k}$.
    \end{lemma}
    \begin{proof}.~
        By assumption, $\ZBold_k$ has full row rank and thus has a right multiplicative inverse $\ZBold_k\pinv$. Hence $\Y_k = \ZBold \ZBold_k\pinv$ satisfies the two relations
        \begin{equation*}
            \U_k\tra \Y_k = \I_k \quad \txtand \quad \overline{\U}_k\tra \Y_k = \T_k.
        \end{equation*}
        Moreover, we have $\range(\Y_k) \subset \range(\ZBold)$ which implies that $\I_n - \Proj(\ZBold) \lleq \I_n - \Proj(\Y_k)$ according to~\cite[Proposition 8.5]{HalkoMartinssonEtAl_FindingStructureRandomness_2011}. By using the conjugation rule~\eqref{eq:conjugation_rule} and the identity $\U_k \U_k\tra + \overline{\U}_k\overline{\U}_k\tra = \I_n$, we obtain
        \begin{align*}
            \U_k\tra [\I_n - \Proj(\ZBold)] \U_k \lleq \U_k\tra [\I_n - \Proj(\Y_k)] \U_k & = \U_k\tra \left(\I_n - \Y_k (\Y_k\tra \Y_k)\inv \Y_k\tra \right) \U_k, \\
            & = \I_k - \U_k\tra \Y_k (\Y_k\tra \Y_k)\inv \Y_k\tra \U_k, \\
            & = \I_k - (\Y_k\tra \Y_k)\inv, \\
            & = \I_k - \left(\Y_k\tra (\U_k \U_k\tra + \overline{\U}_k \overline{\U}_k\tra)\Y_k\right)\inv, \\
            & = \I_k - \left(\I_k + \T_k\tra \T_k \right)\inv, \\
            & = \T_k\tra \left(\I_{n-k} + \T_k \T_k\tra \right)\inv \T_k,
        \end{align*}
        where the last equality is obtained using the Sherman-Morrison formula~\eqref{eq:sherman}. To conclude, we observe that $(\I_{n-k} + \T_k \T_k\tra )\inv \lleq \I_{n-k}$, which implies, using the conjugation rule \eqref{eq:conjugation_rule}, that $\T_k\tra (\I_{n-k} + \T_k \T_k\tra )\inv \T_k \lleq \T_k\tra \T_k$.
    \end{proof}

    \begin{remark}
        Let $\theta_1, \dots, \theta_k$ denote the principal angles between $\range(\U_k)$ and $\range(\ZBold)$, and $\widetilde{\theta}_1, \dots, \widetilde{\theta}_k$ the ones between $\range(\U_k)$ and $\range(\Y_k)$. The statement of Lemma~\ref{lem:1} can be geometrically rephrased as
        \begin{equation*}
            \sin(\theta_i)^2 \leq \sin(\widetilde{\theta}_i)^2 =
            \frac{\tan(\widetilde{\theta}_i)^2}{1 + \tan(\widetilde{\theta}_i)^2} \leq \tan(\widetilde{\theta}_i)^2,
            \quad 1 \leq i \leq k.
        \end{equation*}
        Less formally, the forthcoming analysis is based on how close $\range(\U_k)$ and $\range(\Y_k)$ are, while the truly computed error rather concerns how close $\range(\ZBold)$ is from $\range(\U_k)$. Since the latter is of dimension $\ell$, while $\range(\Y_k)$ is of dimension $k$, the looseness of the derived bounds will increase with $\ell$.
    \end{remark}

    The next result is a deterministic bound for the low-rank approximation error which generalizes~\cite[Proposition 9.6]{HalkoMartinssonEtAl_FindingStructureRandomness_2011}. Indeed, in the particular case where $\ZBold = \A \OmegaBold$ and up to different notation choices, Proposition~\ref{prop:det_analysis} generalizes~\cite[Proposition 9.6]{HalkoMartinssonEtAl_FindingStructureRandomness_2011}.

    \begin{proposition} \label{prop:det_analysis}
        Let $\A \in \Rmn$ be an arbitrary matrix and $\ZBold \in \Rml$ a full column rank matrix such that $\ell \leq \rank(\A)$. Let us further denote $\ZBold_k = \U_k\tra \ZBold \in \R^{k \times \ell}$ and $\overline{\ZBold}_k = \overline{\U}_k\tra \ZBold \in \R^{(m-k) \times \ell}$.

        \vspace{0.25cm}

        \noindent \noindent For any integer $k \leq \ell$, if $\rank(\ZBold_k) = k$, then one has
        \begin{equation*}
            \norm{\left[\I_n - \Proj(\ZBold)\right] \A}_F^2 \leq
            \norm{\overline{\SigmaBold}_k}_F^2 + \norm{\overline{\ZBold}_k \ZBold_k\pinv \SigmaBold_k}_F^2.
        \end{equation*}
    \end{proposition}
    \begin{proof}.~
        From the partitioning of the SVD of $\A$ recalled above, one can define $\A_k = \U_k \SigmaBold_k \V_k\tra$ and $\overline{\A}_k = \overline{\U}_k \overline{\SigmaBold}_k \overline{\V}_k\tra$, from which one gets $\A \A\tra = \A_k \A_k\tra + \overline{\A}_k\overline{\A}_k\tra$. Then, by definition of the Frobenius norm, one has
        \begin{align*}
            \norm{[\I_n - \Proj(\ZBold)] \A}_F^2 & = \tr\left( [\I_n - \Proj(\ZBold)] \A \A\tra [\I_n - \Proj(\ZBold)]\right), \\
            & = \tr\left( [\I_n - \Proj(\ZBold)] [\A_k \A_k\tra + \overline{\A}_k\overline{\A}_k\tra] [\I_n - \Proj(\ZBold)] \right), \\
            & = \tr\left( [\I_n - \Proj(\ZBold)] \A_k \A_k\tra [\I_n - \Proj(\ZBold)] \right) + \tr\left( [\I_n - \Proj(\ZBold)]\overline{\A}_k\overline{\A}_k\tra [\I_n - \Proj(\ZBold)] \right) \\
            & = \norm{[\I_n - \Proj(\ZBold)] \A_k}_F^2 + \norm{[\I_n - \Proj(\ZBold)] \overline{\A}_k}_F^2.
        \end{align*}
        The submultiplicativity of the Frobenius norm \eqref{eq:submultiplicativity} yields $\norm{[\I_n - \Proj(\ZBold)] \overline{\A}_k}_F^2 \leq \norm{[\I_n - \Proj(\ZBold)]}_2^2 \norm{\overline{\A}_k}_F^2$. Then, $\I_n - \Proj(\ZBold)$ is unitary which implies that $\norm{[\I_n - \Proj(\ZBold)]}_2 = 1$ and the unitary invariance of the Frobenius norm yields $\norm{\overline{\A}_k}_F = \norm{\overline{\U}_k \overline{\SigmaBold}_k \overline{\V}_k\tra}_F = \norm{\overline{\SigmaBold}_k}_F$. Overall one gets $\norm{[\I_n - \Proj(\ZBold)] \overline{\A}_k}_F^2 \leq \norm{\overline{\SigmaBold}_k}_F^2$. For the remaining term, one has
        \begin{equation*}
            \norm{[\I_n - \Proj(\ZBold)] \A_k}_F^2 = \tr\left( \SigmaBold_k \U_k\tra [\I_n - \Proj(\ZBold)] \U_k \SigmaBold_k \right).
        \end{equation*}
        Applying Lemma~\ref{lem:1} and conjugating by $\SigmaBold_k$, from~\eqref{eq:conjugation_rule} we get
        \begin{equation*}
            \SigmaBold_k \U_k\tra [\I_n - \Proj(\ZBold)] \U_k \SigmaBold_k \lleq \SigmaBold_k \T_k\tra \T_k \SigmaBold_k.
        \end{equation*}
        Using the monotonicity of the trace~\eqref{eq:monotonic} ends the proof.
    \end{proof}

    If ones consider a low rank approximation of $\A$ of the form $\Proj(\ZBold)\A$, Proposition \ref{prop:det_analysis} shows that the deviation from the optimal low rank approximation error depends on a particular quantity, namely $\norm{\overline{\ZBold}_k \ZBold_k\pinv \SigmaBold_k}_F^2$. Proposition~\ref{prop:stoch} provides a stochastic analysis of this particular quantity, and constitutes the final and main result from which Theorem~\ref{th:main:general} can be deduced.
    
    \begin{proposition} \label{prop:stoch}
        Let $\A \in \Rmn$ be an arbitrary matrix, and $\ZBold = \D \G \in \Rml$ where $\D \in \Rmr$ and $\G \in \Rrl$ is a matrix whose columns are independently sampled from an $r$-variate standard Gaussian distribution, and such that $\ell \leq \min(\rank(\A), \, \rank(\D))$. Let us further consider the matrices $\U_k$, $\overline{\U}_k$, $\SigmaBold_k$ and $\overline{\SigmaBold}_k$ related to the singular value decomposition of $\A$ and defined in \eqref{eq:partitioning_svd}, and the matrices $\ZBold_k = \U_k\tra \ZBold \in \R^{k \times \ell}$ and $\overline{\ZBold}_k = \overline{\U}_k\tra \ZBold \in \R^{(m-k) \times \ell}$.
        
        \vspace{0.25cm}

        \noindent For any integer $\ell \geq k + 2$, if $\rank(\D\tra \U_k) = k$, then one has
        
        \begin{equation*}
            \expect{\norm{\overline{\ZBold}_k \ZBold_k\pinv \SigmaBold_k}_F^2} =
            \left(\tau_k^2 + \dfrac{\rho_k^2}{\ell - k + 1}\right) \norm{\overline{\SigmaBold}_k}_F^2 .
        \end{equation*}

        \noindent Moreover, if $\ell \geq k + 4$, then for all $u, t \geq 1$ one has
        \begin{equation*}
            \norm{\overline{\ZBold}_k \ZBold_k\pinv \SigmaBold_k}_F \leq
            \left(\tau_k + \dfrac{\rho_k}{\ell - k - 1}\right)
            \norm{\overline{\SigmaBold}_k}_F,
        \end{equation*}
        holds with probability at least $1 - e^{-u^2 / 2} - t^{-(\ell - k)}$, with 
        \begin{align*}
            & \tau_k = \frac{1}{\norm{\overline{\SigmaBold}_k}_F} \,
            \norm{\tanmat(\U_k, \, \D\D\tra \U_k) \SigmaBold_k}_F, \\ \txtand \quad
            & \rho_k = \frac{1}{\norm{\overline{\SigmaBold}_k}_F} \,
            \norm{[\I_r - \Proj(\D\tra \U_k)] \D\tra \overline{\U}_k}_F \,
            \norm{\SigmaBold_k (\U_k\D\D\tra\U_k)\halfinv}_F.
        \end{align*}
    \end{proposition}
    \begin{proof}.
        ~
        Let us define $\K = \D\D\tra \in \Rmm$ so that $\K$ is the covariance matrix of the columns of $\ZBold$.
        \paragraph*{Expectation.}
            Applying the theorem of total expectation yields
            \begin{equation*}
                \expect{\norm{\overline{\ZBold}_k \ZBold_k\pinv \SigmaBold_k}_F^2} =
                \expect{\expect{\norm{\overline{\ZBold}_k \ZBold_k\pinv \SigmaBold_k}_F^2 \mid \ZBold_k}}.
            \end{equation*}
            Let us consider the following partitioning
            \begin{equation*}
                \U\tra \K \U = \left[\U_k \mid \overline{\U}_k\right]\tra \K \left[\U_k \mid \overline{\U}_k\right] =
                \begin{bmatrix}
                    \K_k & \K_{\perp,k}\tra \, \\
                    \,\K_{\perp,k} & \overline{\K}_k
                \end{bmatrix}.
            \end{equation*}
            In the inner expectation, $\ZBold_k\pinv$ is a fixed matrix, but unlike the standard result in~\cite{HalkoMartinssonEtAl_FindingStructureRandomness_2011}, $\ZBold_k$ and $\overline{\ZBold}_k$ are not statistically independent here. This implies that conditioned by $\ZBold_k$, the columns of $\overline{\ZBold}_k$ no longer follow a standard multivariate Gaussian distribution. Instead, if $[\ZBold_k]_{(i)}$ denotes the $i$-th column of $\ZBold_k$, then the $i$-th column of $\overline{\ZBold}_k$ follows a $(m-k)$-variate Gaussian distribution with mean term
            \begin{equation*}
                \bm{\mu}_i \equiv \K_{\perp,k} \K_k\inv [\ZBold_k]_{(i)},
            \end{equation*}
            and covariance matrix
            \begin{equation*}
                \K / \K_k \equiv  \overline{\K}_k - \K_{\perp,k} \K_k\inv \K_{\perp,k}\tra.
            \end{equation*}
            Consequently, if one defines $\bm{\mu} = [\bm{\mu}_1, \, \dots, \, \bm{\mu}_\ell] \in \R^{(m-k) \times \ell}$, then the distribution of $\overline{\ZBold}_k$ conditioned by $\ZBold_k$ is the same as $\bm{\mu} + (\K / \K_k)\half \G$ where the columns of $\G \in \R^{(m-k) \times \ell}$ are sampled from a standard Gaussian distribution. Using this fact, one has
            \begin{align*}
                \norm{[\bm{\mu} + (\K / \K_k)\half \G] \ZBold_k\pinv \SigmaBold_k}_F^2 & = \norm{\bm{\mu} \, \ZBold_k\pinv \SigmaBold_k + (\K / \K_k)\half \G \ZBold_k\pinv \SigmaBold_k}_F^2 \\
                & = \norm{\bm{\mu} \, \ZBold_k\pinv \SigmaBold_k}_F^2 + \norm{(\K / \K_k)\half \G \ZBold_k\pinv \SigmaBold_k}_F^2 +
                2\tr\left(\SigmaBold_k [\ZBold_k\pinv]\tra \bm{\mu}\tra (\K / \K_k)\half \G \ZBold_k\pinv \SigmaBold_k \right).
            \end{align*}
            The linearity of the expectation allows us to deal with each term separately. Conditioned by $\ZBold_k$, the first one is a non-random constant, and the second term can be handled using~\cite[Proposition 10.1]{HalkoMartinssonEtAl_FindingStructureRandomness_2011} since $\G$ has a standard multivariate Gaussian distribution. Concerning the third term, the linearity of the expectation conditioned by $\ZBold_k$ yields
            \begin{equation*}
            \expect{\tr\left(\SigmaBold_k [\ZBold_k\pinv]\tra \bm{\mu}\tra (\K / \K_k)\half \G \ZBold_k\pinv \SigmaBold_k \right) \mid \ZBold_k} =
                \tr\left(\SigmaBold_k [\ZBold_k\pinv]\tra \bm{\mu}\tra (\K / \K_k)\half \expect{\G \mid \ZBold_k}\ZBold_k\pinv \SigmaBold_k \right)= 0.
            \end{equation*}
            Altogether, one has
            \begin{align*}
                \expect{\norm{\overline{\ZBold}_k \ZBold_k\pinv \SigmaBold_k}_F^2 \mid \ZBold_k} & =
                \norm{\bm{\mu} \, \ZBold_k\pinv \SigmaBold_k}_F^2 +
                \norm{(\K / \K_k)\half}_F^2 \, \norm{\ZBold_k\pinv \SigmaBold_k}_F^2 \\
                & = \norm{\K_{\perp,k} \K_k\inv \SigmaBold_k}_F^2 +
                \norm{(\K / \K_k)\half}_F^2 \, \norm{\ZBold_k\pinv \SigmaBold_k}_F^2,
            \end{align*}
            where the last equality follows from the simplification of $\ZBold_k \ZBold_k\pinv$ when substituting $\bm{\mu}$ by its expression. Taking again the expectation yields
            \begin{equation*}
                \expect{\norm{\overline{\ZBold}_k \ZBold_k\pinv \SigmaBold_k}_F^2} =
                \norm{\K_{\perp,k} \K_k\inv \SigmaBold_k}_F^2 +
                \norm{(\K / \K_k)\half}_F^2 \cdot \expect{\norm{\ZBold_k\pinv \SigmaBold_k}_F^2}.
            \end{equation*}   
            The columns of $\ZBold_k$ are independently sampled from a $k$-variate centered Gaussian distribution with covariance matrix $\K_k$, with $\K_k = \U_k\tra \D\D\tra \U_k$ being of rank $k$ by assumption.
            Consequently, the matrix $\ZBold_k \ZBold_k\tra$ follows a Wishart distribution with $\ell-1$ degrees of freedom and covariance matrix $\K_k$~\cite[Definition 3.1.3]{Muirhead_AspectsMultivariateStatistical_1982}.
            Since $\ell \geq k + 1$ by assumption, $\ZBold_k \ZBold_k\tra$ is non-singular with probability 1~\cite[Theorem 3.1.4]{Muirhead_AspectsMultivariateStatistical_1982}, implying that $\ZBold_k\pinv = \ZBold_k\tra(\ZBold_k\ZBold_k\tra)\inv$.
            Thus, one has $\norm{\ZBold_k\pinv \SigmaBold_k}_F^2 = \tr(\SigmaBold_k [\ZBold_k\pinv]\tra \ZBold_k\pinv \SigmaBold_k) = \tr(\SigmaBold_k [\ZBold_k\ZBold_k\tra]\inv \SigmaBold_k)$. Then, as a consequences of~\cite[Theorem 3.2.12]{Muirhead_AspectsMultivariateStatistical_1982}, since $\ell \geq k + 2$ by assumption, one has
            \begin{equation*}
                \expect{[\ZBold_k\ZBold_k\tra]\inv} =
                \dfrac{\SigmaBold_k \K_k \inv \SigmaBold_k}{\ell - k - 1}.
            \end{equation*}
            
            Using the linearity of the expectation finally yields
            \begin{equation*}
                \expect{\norm{\ZBold_k\pinv \SigmaBold_k}_F^2} =
                \expect{\tr(\SigmaBold_k [\ZBold_k \ZBold_k\tra]\inv \SigmaBold_k)} =
                \tr(\SigmaBold_k \expect{[\ZBold_k \ZBold_k\tra]\inv} \SigmaBold_k) =
                \dfrac{\tr(\SigmaBold_k \K_k\inv \SigmaBold_k)}{\ell - k - 1} =
                \dfrac{\norm{\SigmaBold_k \K_k\halfinv}_F^2}{\ell - k - 1}.
            \end{equation*}
            Altogether, one obtains,
            \begin{equation*}
                \expect{\norm{\overline{\ZBold}_k \ZBold_k\pinv \SigmaBold_k}_F^2} =
                \norm{\K_{\perp,k} \K_k\inv \SigmaBold_k}_F^2 +
                \norm{(\K / \K_k)\half}_F^2 \cdot \dfrac{\norm{\SigmaBold_k \K_k\halfinv}_F^2}{\ell - k - 1}.
            \end{equation*}  
            Up to a rearrangement of $\K_{\perp,k} \K_k\inv$ and $\norm{(\K / \K_k)\half}_F^2$ proved below, this ends the part of the proof related to the expectation.

        \paragraph*{Probability.}
            Let $\varphi : \R^{(m-k) \times \ell} \rightarrow \R, \, \J \mapsto \norm{\, [\bm{\mu} + (\K / \K_k)\half \J] \, \ZBold_k\pinv \SigmaBold_k}_F$. Using the reverse triangular inequality, it holds for any $\J_1, \J_2 \in \R^{(m-k) \times \ell}$ that
            \begin{align*}
                \abs{\varphi(\J_1) - \varphi(\J_2)} & \leq
                \norm{(\K / \K_k)\half \, [\J_1 - \J_2] \, \ZBold_k\pinv \SigmaBold_k}_F \\
                & \leq
                \norm{(\K / \K_k)\half}_F \norm{\ZBold_k\pinv \SigmaBold_k}_F \, \norm{[\J_1 - \J_2]}_F,
            \end{align*}
            where the last inequality follows from the submultiplicativity property. Thus, $\varphi$ is at worst $L$-Lipschitz with $L = \norm{(\K / \K_k)\half}_F \norm{\ZBold_k\pinv \SigmaBold_k}_F$. Applying~\cite[Proposition 10.3]{HalkoMartinssonEtAl_FindingStructureRandomness_2011} then yields that if $\G \in \R^{(m-k) \times \ell}$ is a matrix whose columns are sampled from a standard Gaussian distribution, then $\forall \, u \geq 0$ one has
            \begin{equation*}
                \P\set{\varphi(\G) \leq \expect{\varphi(\G)} + Lu} \geq 1 - e^{-u^2 / 2}.
            \end{equation*}
            Hölder's inequality yields $\expect{\varphi(\G)} \leq \mathbb{E}[\varphi(\G)^2]\half$, and combining this fact with the result proved right above and the properties of multivariate Gaussian distribution one has
            \begin{align*}
                \expect{\varphi(\G)} & \leq \left(\expect{\norm{\overline{\ZBold}_k \ZBold_k\pinv \SigmaBold_k}_F^2 \mid \ZBold_k}\right)\half \\ 
                & = \left(\norm{\K_{\perp,k} \K_k\inv \SigmaBold_k}_F^2 + \norm{(\K / \K_k)\half}_F^2 \, \norm{\ZBold_k\pinv \SigmaBold_k}_F^2\right)\half \\
                & \leq \norm{\K_{\perp,k} \K_k\inv \SigmaBold_k}_F + \norm{(\K / \K_k)\half}_F \, \norm{\ZBold_k\pinv \SigmaBold_k}_F.
            \end{align*}
            Altogether, this implies that
            \begin{equation} \label{eq:EZk}
                \P\set{\norm{\overline{\ZBold}_k \ZBold_k\pinv \SigmaBold_k}_F \leq
                \norm{\K_{\perp,k} \K_k\inv \SigmaBold_k}_F + (1+u) \cdot \norm{(\K / \K_k)\half}_F \norm{\ZBold_k\pinv \SigmaBold_k}_F} \geq 1 - e^{-u^2 / 2}.
            \end{equation}
            Let us now consider, for $t \geq 1$, the event
            \begin{equation*}
                E_t = \set{\norm{\ZBold_k\pinv \SigmaBold_k}_F \leq \sqrt{3} t \cdot \sqrt{\frac{\tr(\SigmaBold_k^2 \K_k\inv)}{\ell - k + 1}}}.
            \end{equation*}
            Applying~\cite[Lemma 8]{BoulleTownsend_GeneralizationRandomizedSingular_2022} ($\ell \geq k + 4$) yields for all $t \geq 1$ that $\P\set{E_t^C} \leq t^{-(\ell - k)}$ where $E_t^C$ denotes the complement of the event set $E_t$. Now, by denoting $E_{\ZBold_k}$ the event considered in~\eqref{eq:EZk}, the law of total probabilities reads
            \begin{equation*}
                \P\set{E_{\ZBold_k}} = \P\set{E_{\ZBold_k} \mid E_t} \P\set{E_t} + \P\set{E_{\ZBold_k} \mid E_t^{\mathsf{C}}} \P\set{E_t^{\mathsf{C}}},
            \end{equation*}
            Since one trivially has $\P\set{E_t} \leq 1$ and $\P\set{E_{\ZBold_k} \mid E_t^{\mathsf{C}}} \leq 1$, one gets
            \begin{equation*}
                1 - e^{-u^2 / 2} \leq \P\set{E_{\ZBold_k}} \leq \P\set{E_{\ZBold_k} \mid E_t} + t^{-(\ell - k)},
            \end{equation*}
            or equivalently,
            \begin{equation*}
                \P\set{E_{\ZBold_k} \mid E_t} \geq 1 - e^{-u^2 / 2} - t^{-(\ell - k)},
            \end{equation*}
            which ends the proof for the bound in probability.

        \paragraph*{Rearranging the terms.}
            It remains to transform the terms to get more interpretable bounds. Since by definition $\left[\U_k \mid \overline{\U}_k\right]$ is orthogonal, and $\ZBold$ is full column rank, then using~\eqref{eq:tangent}, one readily has
            \begin{equation*}
                \K_{\perp,k} \K_k\inv =
                \overline{\U}_k\tra \K \U_k (\U_k\tra \K \U_k)\inv =
                \tanmat(\U_k, \K \U_k) =
                \tanmat(\U_k, \D\D\tra \U_k).
            \end{equation*}
            Then, using classical algebraic manipulations, one finally gets
            \begin{align*}
                \norm{(\K / \K_k)\half}_F^2 & =
                \tr\left(\overline{\U}_k\tra \K \overline{\U}_k - \overline{\U}_k\tra \K \U_k \K_k\inv \U_k\tra \K \overline{\U}_k\right) \\
                & = \tr\left(\overline{\U}_k\tra \D \left[\I_r - \D\tra \U_k \K_k\inv \U_k\tra \D \right] \D\tra\overline{\U}_k\right) \\
                & = \tr\left(\overline{\U}_k\tra \D \left[\I_r - \Proj(\D\tra \U_k) \right] \D\tra\overline{\U}_k\right)
                = \norm{[\I_r - \Proj(\D\tra \U_k)] \D\tra\overline{\U}_k}_F^2,
            \end{align*}  
            which concludes the proof.
        \end{proof}

    Theorem~\ref{th:main:general} can now be proved as an almost straightforward combination of Propositions \ref{prop:det_analysis} and \ref{prop:stoch}.

    \paragraph*{Proof of Theorem~\ref{th:main:general}}
        The expectation is monotonic and linear, so taking the expectation over $\ZBold$ in Proposition~\ref{prop:det_analysis} yields
        \begin{equation*}
            \expect{\norm{\left[\I_n - \Proj(\ZBold)\right] \A}_F^2} \leq
            \norm{\overline{\SigmaBold}_k}_F^2 + \expect{\norm{\overline{\ZBold}_k \ZBold_k\pinv \SigmaBold_k}_F^2}.
        \end{equation*}
        Then, since $\ZBold = \D \G \in \Rml$ satisfies the assumptions of Proposition \ref{prop:stoch}, applying the result for the expectation yields the desired bound.
        \begin{equation*}
            \expect{\norm{\left[\I_n - \Proj(\ZBold)\right] \A}_F^2} \leq
            \left(1 + \tau_k^2 + \frac{\rho_k^2}{\ell - k - 1}\right) \norm{\overline{\SigmaBold}_k}_F^2.
        \end{equation*}
        Using Hölder's inequality, that is $\mathbb{E}[\,\, \abs{X}^2] \leq \mathbb{E}[\, \abs{X}]^2$ for any random variable $X$, concludes the proof. For the probability bounds, using the fact that $\sqrt{a + b} \leq \sqrt{a} + \sqrt{b}$ for any $a, b \geq 0$, one can deduce from Proposition~\ref{prop:det_analysis} that
        \begin{equation*}
            \norm{\left[\I_n - \Proj(\ZBold)\right] \A}_F \leq \norm{\overline{\SigmaBold}_k}_F + \norm{\overline{\ZBold}_k \ZBold_k\pinv \SigmaBold_k}_F.
        \end{equation*}
        Then, applying the result in probability of Proposition~\ref{prop:stoch} on the second term of the right-hand side concludes the proof.
\bibliographystyle{etna}
\bibliography{bibliography.bib}

\end{document}